%% file: main.tex
\documentclass[final,a4paper,leqno]{siamltex}

\usepackage[english]{babel}
\usepackage{amsmath,amssymb,amsfonts,graphicx,subfigure}
\usepackage[ruled]{algorithm2e} 
\usepackage{myshortcuts}
\newtheorem{remark}[theorem]{Remark} 
\usepackage{mathabx}
\usepackage{url}

\usepackage{pgfplots}
\usetikzlibrary{external,positioning, fit}
\usetikzlibrary{decorations.pathreplacing,patterns}

\iftrue 
  \usepackage{tikz}
  \usepackage{pgfplots}
  \pgfplotsset{
    compat=newest,
    tick label style={font=\scriptsize},
    label style={font=\scriptsize},
    legend style={font=\scriptsize}
  }
  \usetikzlibrary{decorations.pathreplacing}
  \iffalse 
     \usepgfplotslibrary{external} 
  \else
     \renewcommand{\tikzsetnextfilename}[1]{}
  \fi
\else
  \usepackage{tikzexternal}
  \tikzsetexternalprefix{gfxtmp/}
  
\fi

\newcommand{\coloneqq}{\mathrel{\mathop:}=}
\newcommand{\pideg}[1]{\Pi_{#1}}
\newcommand{\pigrad}[1]{\widebar{\Pi}_{#1}}
\newcommand{\pistar}[1]{\Pi_{2^{#1}}^*}
\newcommand{\pistarbar}[1]{\widebar{\Pi_{2^{#1}}^*}}

\usepackage{blkarray, makecell}
\usepackage{booktabs}

\title{The polynomial set associated with a fixed number of matrix-matrix multiplications}

\author{Elias Jarlebring and Gustaf Lorentzon}
\date{\today}

\begin{document}
\maketitle
\begin{abstract}
We consider the problem of computing matrix polynomials $p(X)$, where $X$ is a large dense matrix,
with as few matrix-matrix multiplications as possible. More precisely, let $\pistar{m}$ represent the set of
polynomials computable with $m$ matrix-matrix multiplications, but with an arbitrary number of  matrix additions and scaling operations.
We characterize this set through a tabular parameterization.
By deriving equivalence transformations of the tabular representation,
we establish new methods that can be used to construct elements 
of $\pistar{m}$
and determine general properties of the set.
The transformations allow us to eliminate variables and prove that the dimension is bounded by $m^2$, which is subsequently proven to be sharp, i.e., $\dim(\pistar{m})=m^2$.
Consequently, we have identified a parameterization that, to the best of our knowledge,
is the first minimal parameterization. We also conduct a study using computational tools from algebraic geometry
to determine the largest degree $d$ such that all polynomials of that degree
belong to $\pistar{m}$, or its closure. In many cases, the computational setup is constructive in the sense that it can also be used to determine a specific evaluation scheme for a given polynomial.

\end{abstract}
\section{Introduction}
\input{intro}

\section{Evaluation scheme}\label{sec:eval}

\subsection{Definition of the evaluation scheme}

Fundamental to the results of this paper is an evaluation scheme for constructing elements of $\pistar{m}$. Such evaluation schemes have been presented in various forms in \cite[Section~2]{Paterson:1973:OnTheNumber}, \cite[section~3]{bbc19}, and \cite{jarlebring:2023:graphs}, where in \cite{jarlebring:2023:graphs} they are referred to as degree-optimal polynomials. The evaluation scheme involves parameters stored in the matrices $A$ and $B$, and the 
vector $c$. Matrices $A$ and $B$ contain coefficients for linear combinations, where each row $k$ corresponds to the linear combinations associated with the $k$th multiplication. As we perform $m$ multiplications, these matrices have $m$ rows. The vector $c$ contains coefficients for linear combinations that are used after the multiplications have been completed.

In particular, consider the triplet
$(A,B,c)\in\CC^{m\times (m+1)}\times\CC^{m\times (m+1)}\times\CC^{m+2}$
and associate the following sequence of operations involving exactly $m$ matrix-matrix multiplications. Define $Q_1 = I$ and $Q_2 = X$. Then, we iterate the process to generate $Q_3,\ldots,Q_{m+2}$ using the elements of matrices $A$ and $B$ as follows:
\begin{equation}
\begin{alignedat}{3}\label{eq:generateB}
Q_3 &=(a_{11}Q_1 +a_{12}Q_2 )(b_{11}Q_1 +b_{12}Q_2 ) \\
Q_4 &=(a_{21}Q_1 +a_{22}Q_2 +a_{23}Q_3 )(b_{21}Q_1 +b_{22}Q_2 +b_{23}Q_3 ) \\
Q_5 &=(a_{31}Q_1 +a_{32}Q_2 +a_{33}Q_3 +a_{34}Q_4 )(b_{31}Q_1 +b_{32}Q_2 +b_{33}Q_3 +b_{34}Q_4 ) \\ 
& \vdots \\
Q_{m+2} &=(a_{m,1}Q_1 + \cdots + a_{m,m+1} Q_{m+1} )(b_{m,1}Q_1 + \cdots + b_{m,m+1} Q_{m+1} ).
\end{alignedat}        
\end{equation}
Furthermore, we compute the output of the scheme by forming the linear combination corresponding to the $c$-vector:
\begin{equation}\label{eq:poly_def}
p(X) = c_1 Q_1 + c_2 Q_2 + \cdots + c_{m+2} Q_{m+2}.
\end{equation}
Thus, a given triple $(A,B,c)$, with $A$ and $B$ being matrices, defines a polynomial $p\in\pistar{m}$.

Note that all evaluation schemes can be expressed in this manner, due to the fact that for each multiplication step we use all preceding information available. For instance, in the second multiplication step we determine $Q_4$ by computing the product of two linear combinations of $Q_3$, $Q_2$ and $Q_1$. Hence, a triplet $(A, B, c)$ parameterizes $\pistar{m}$ through the expressions in \eqref{eq:generateB} and \eqref{eq:poly_def}.

\subsection*{Examples of instances of $\pistar{m}$}
To see how the tables represent standard evaluation schemes, such as monomial evaluation, Horner evaluation, and Paterson--Stockmeyer evaluation, see \cite[Section~3]{jarlebring:2023:graphs}.
In the following, we show how to evaluate 
\[
p(X)=X^7 + \epsilon X^8,
\]
in three multiplications. 
Consider the coefficient triplet:
\begin{subequations} \label{eq:epsilonexample}
\begin{align}
        [A|B] &= \left[\begin{array}{c c c c | c c c c}
            0 & 1 &   &   & 0 & 1 &   &  \\
            0 & \frac{1}{2\epsilon} & 1 &   & 0 & 0 & 1 &  \\
            0 & \frac{7 + 8 \epsilon^2 -16\epsilon^4}{128\epsilon^5} & -\frac{1}{8\epsilon^2} - \frac{1}{2} & 1 & 0 & \frac{-7 + 8 \epsilon^2 + 16 \epsilon^4}{128\epsilon^5} & -\frac{1}{8 \epsilon^2} + \frac{1}{2} & 1
            \end{array}\right] \label{eq:epsilonexampleAB}\\
            c &= \left[ \begin{array}{ccccc}
                 0 & 0 & \frac{49 - 288 \epsilon^4 + 256 \epsilon^8 }{16384 \epsilon^9} & \frac{-5 + 16 \epsilon^4}{ 64 \epsilon^3 } & \epsilon
            \end{array}\right]. \label{eq:epsilonexampleC}
\end{align}
\end{subequations}

For this triplet $(A,B,c)$, we obtain the polynomial $p(X) = X^7 + \epsilon X^8$
using the evaluations specified in 
\eqref{eq:generateB}-\eqref{eq:poly_def}.
We observe that $X^7$ is a limit point of $\pistar{3}$ corresponding to
$\epsilon\rightarrow 0$,
even though $X^7$ requires four multiplications.
This example illustrates the complexity of $\pistar{m}$ for $m=3$.
As $m$ increases, the complexity becomes even more intricate.

 \input{variety.tex}
\medskip\medskip\medskip
\section{Equivalence transformations}\label{sec:transform}
\subsection{Substitution transformations}

The following theorems describe transformations of the evaluation 
scheme $(A,B,c)$, as defined in Section~\ref{sec:eval},
such that the output polynomial $p(X)$ is unchanged. 
The modified evaluation scheme will be denoted $(\hat{A},\hat{B},\hat{c})$ and the corresponding $Q$-matrices will be denoted 
$\hat{Q}_1,\ldots,\hat{Q}_{m+1}$.
The first transformation can be viewed as a change of variables, where we rescale one 
of the $Q$-coefficients. For example if we set
\[
\hat{Q}_4=\alpha Q_4
\]
the tables \eqref{eq:generateB} are changed by
\[
\hat{Q}_4=(\alpha a_{21}Q_1 +\alpha a_{22}Q_2 +\alpha a_{23}Q_3            )( b_{21}Q_1  +b_{22}Q_2 +b_{23}Q_3).
\]
In order to keep $Q_5,\ldots, Q_{m+2}$ unchanged, we can cancel the transformation by scaling the coefficient in the corresponding term, e.g.,
\[
Q_5=
(a_{31}Q_1 + a_{32}Q_2 + a_{33}Q_3+\alpha^{-1}a_{34}\hat{Q}_4)
(b_{31}Q_1 + b_{32}Q_2 + b_{33}Q_3+\alpha^{-1}b_{34}\hat{Q}_4).
\]
This implies that without changing the polynomial $p$, we may rescale one row in $A$ and one column in $A$ and $B$. This is formalized in the following
theorem, where we omit indices on elements of $A$ and $B$ that are unchanged in
$\hat{A}$ and $\hat{B}$. By symmetry, we note that the same result holds for $A$ and $B$ switched.

\begin{theorem}\label{thm:rowscaling}
    Let $p$ be the polynomial associated with $(A,B,c)\in \CC^{m \times (m+1)}\times\CC^{m \times (m+1)}\times \CC^{m+2}$ and let $\hat{p}$ be the polynomial associated with $(\hat{A},\hat{B},\hat{c})$, where
        \begin{equation}\label{eq:rowscaling_ABc}
            \begin{alignedat}{4}
        &\hat{A} &&=
            \left[\begin{array}{c c c c c c c c c c}
                a               & a      & \\
                \vdots          &        & \ddots \\
                a               & \cdots & \cdots   & a    \\
                \hat{a}_{k,1}   & \cdots & \cdots   & \cdots    & \hat{a}_{k,k+1} \\
                a               & \cdots & \cdots   & \cdots    & a     & \hat{a}_{k+1,k+2}  \\
                \vdots          &        &          &           & \vdots& \vdots & a \\
                \vdots          &        &          &           & \vdots& \vdots & \vdots &\ddots \\
                a               & \cdots &  \cdots  & \cdots    & a     & \hat{a}_{m,k+2} & a  &\cdots  & a \\
            \end{array}\right] \\
        &\hat{B} &&=
            \left[\begin{array}{c c c c c c c c c c}
                b               & b      & \\
                \vdots          &        & \ddots \\
                \phantom{x}b\phantom{x}               & \cdots & \cdots   & b    \\
                b   & \cdots & \cdots   & \cdots    & \phantom{xx}b\phantom{xx} \\
                b               & \cdots&\cdots   & \cdots    & b     & \hat{b}_{k+1,k+2}  \\
                \vdots          & &         &           & \vdots& \vdots & b \\
                \vdots          & &         &           & \vdots& \vdots & \vdots &\ddots \\
                b               & \cdots & \cdots  & \cdots    & b     & \hat{b}_{m,k+2} & b  &\cdots  & b \\
            \end{array}\right] \\
    %
    %
            &\hat{c} &&= 
                \;\left[
                    \begin{array}{c c c c c c c c c}
                        \phantom{x}c_1 & \cdots & \cdots & \cdots & c_{k+1} & \alpha^{-1}c_{k+2} & c_{k+3} &\cdots & c_{m+2} 
                    \end{array}
                \right],
            \end{alignedat}
        \end{equation}
        and
        \begin{subequations}    
        \begin{eqnarray}
            \hat{a}_{k,j}&=&\alpha a_{k,j}, \;\;\; \qquad j=1,\ldots,k+1, \label{eq:rowscaling_akj_update}\\
            \hat{a}_{i,k+2}&=&\alpha^{-1} a_{i,k+2},\;\;\; i=k+1,\ldots,m, \label{eq:rowscaling_akj2_update}\\
            \hat{b}_{i,k+2}&=&\alpha^{-1} b_{i,k+2},\;\;\; i=k+1,\ldots,m. 
        \end{eqnarray}
        \end{subequations}
        Then,
        \begin{equation}\label{eq:scalingproof:p_equal}
            \hat{p} = p,
        \end{equation}
        for any $k\in\left[1,\ldots,m\right]$ and $\alpha \in \CC \backslash \{0\}$. \\      
\end{theorem}
\begin{proof}
By definition, we have $\hat{Q}_1 = Q_1 = I$ and $\hat{Q}_2 = Q_2 = X$. The proof is based
on establishing the following three statements:
\begin{subequations}\label{eq:scalingproof:Bcases}
\begin{eqnarray}
    \hat{Q}_{i+2}   &=& Q_{i+2} \qquad \text{for } i = 1, \ldots, k-1, \label{eq:scalingproof:Bcases1}\\
    \hat{Q}_{k+2}   &=& \alpha Q_{k+2}, \label{eq:scalingproof:Bcases2}\\
    \hat{Q}_{i+2}   &=& Q_{i+2} \qquad \text{for } i = k+1, \ldots, m+2. \label{eq:scalingproof:Bcases3}
\end{eqnarray}
\end{subequations}
Together with the definition of $\hat{c}$ in \eqref{eq:rowscaling_ABc}, this implies that the conclusion of the theorem, stated in \eqref{eq:scalingproof:p_equal}, holds.

To prove \eqref{eq:scalingproof:Bcases1}, we observe that since the first $k-1$ rows of $\hat{A}$ and $\hat{B}$ are unchanged, the variables $\hat{Q}_1,\ldots, \hat{Q}_{k+1}$ are also unchanged.

To prove \eqref{eq:scalingproof:Bcases2}, we consider row $k$, i.e., the first row where there are changes. 
Inserting \eqref{eq:rowscaling_akj_update} in the definition
of $\hat{Q}_{k+1}$ yields a scaling,
    \begin{alignat*}{5}
        \hat{Q}_{k+2} &=( \alpha a_{k,1}Q_1 
        + \cdots + \alpha &&a_{k,k+1} Q_{k+1} &&)(b_{k,1}Q_1 + \cdots + b_{k,k+1} Q_{k+1} &&) =\alpha Q_{k+2},
%
    \end{alignat*}
    which proves \eqref{eq:scalingproof:Bcases2}.
    
    To prove \eqref{eq:scalingproof:Bcases3}, we first consider \eqref{eq:scalingproof:Bcases3} for $i = k+1$.
    Inserting \eqref{eq:rowscaling_akj2_update} in the definition of
    $\hat{Q}_{k+3}$ results in a cancellation 
    \begin{equation*}\label{eq:scalingproof:b>kp2}
        \begin{alignedat}{5}
            \hat{Q}_{k+3} &=(a_{k+1,1}Q_1 + \hdots + \alpha^{-1} a_{k+1,k+2} \hat{Q}_{k+2} &&)(b_{k+1,1}Q_1 + \hdots + \alpha^{-1} b_{k+1,k+2} \hat{Q}_{k+2} &&) \\
              &= Q_{k+3},
        \end{alignedat}
    \end{equation*}
    where the last equality follows from \eqref{eq:scalingproof:Bcases2}. The general statement \eqref{eq:scalingproof:Bcases3} follows from induction.
\end{proof}

In a similar fashion, we carry out a change of variables 
where we add a given value $\alpha$ to one of the elements in the first column. For example, adding $\alpha$ to the 
coefficient corresponding to $Q_1=I$ in the first multiplication leads to
\[
  \hat{Q}_3=
  \left((a_{11}+\alpha)Q_1+a_{12}Q_2\right)
  \left(b_{11}Q_1+b_{12}Q_2\right),
\]
which can be expanded as
\begin{equation}\label{eq:Q3_transf2}
\hat{Q}_3=Q_3+\alpha (b_{11}Q_1+b_{12}Q_2).
    \end{equation}
Let $A_k$ and $B_k$ be the factors that form $Q_k$, i.e., $Q_k=A_kB_k$.
In order to keep $Q_4,\ldots, Q_{m+2}$ unchanged,  
we keep both factors in $\hat{Q}_4=Q_4=A_4B_4$ unchanged by compensating for the transformation \eqref{eq:Q3_transf2}:
\begin{eqnarray}
    A_4&=&(a_{21}-\alpha a_{23}b_{11}) Q_1+
(a_{22}-\alpha a_{23}b_{12}) Q_2+a_{23}\hat{Q}_3\\
    B_4&=&(b_{21}-\alpha b_{23}b_{11}) Q_1+
(b_{22}-\alpha b_{23}b_{12}) Q_2+b_{23}\hat{Q}_3.
\end{eqnarray}

Hence, a modification in the coefficient corresponding to $Q_1$ can be compensated for by modifying the coefficients in all rows below the modification.
This can be applied to an arbitrary row $k$ and arbitrary $\alpha$.

\begin{theorem}\label{thm:1stCol}
Let $p$ be the polynomial associated with $(A,B,c)\in \CC^{m \times (m+1)} \times \CC^{m \times (m+1)} \times \CC^{m+2}$ and let $\hat{p}$ be the polynomial associated with $(\hat{A},\hat{B},\hat{c})$ where
        \begin{equation}\label{eq:1stcol_ABc}
            \begin{alignedat}{4}
        &\hat{A} &&=
            \left[\begin{array}{c c c c c c c c c c}
                a               & a      & \\
                \vdots          &        & \ddots \\
                a               & a & \cdots   & a    \\
                a + \alpha   & a & \cdots   & \cdots    & a \\
                \hat{a}_{k+1,1} & \hat{a}_{k+1,2} & \cdots   & \cdots    & \hat{a}_{k+1,k+1}  & a  \\
                \vdots          & \vdots          &          &           & \vdots& \vdots &\ddots \\
                \hat{a}_{m,1}   & \hat{a}_{m,2}   &  \cdots  & \cdots    & \hat{a}_{m,k+1}  & a  &\cdots  & a \\
            \end{array}\right] \\
        &\hat{B} &&=
            \left[\begin{array}{c c c c c c c c c c}
                b               & b      & \\
                \vdots          &        & \ddots \\
                b               & \cdots & \cdots   & b    \\
                \hat{b}_{k+1,1} & \cdots   & \cdots    & \hat{b}_{k+1,k+1}  & b  \\
                \vdots          &         &           & \vdots& \vdots &\ddots \\
                \hat{b}_{m,1}   & \cdots  & \cdots    & \hat{b}_{m,k+1}  & b  &\cdots  & b \\
            \end{array}\right] \\
    %
    %
            &\hat{c} &&= 
                \;\left[\;
                    \begin{array}{c c c c c c c c}
                        \phantom{x}\hat{c}_1\phantom{x} & \cdots & \cdots &\phantom{x}\hat{c}_{k+1}\phantom{x} & c & \cdots & \cdots & c
                    \end{array}
                \right].
            \end{alignedat}
        \end{equation}
        and
        \begin{subequations}\label{eq:1stColproof:Coeffs}
        \begin{eqnarray}
            \hat{a}_{i,j}&=& a_{i,j} - \alpha a_{i,k+2}b_{k,j},\; \quad j=1,\ldots,k+1, \quad i=k+1,\ldots,m, \label{eq:1stColproof:Coeffs1}\\
            \hat{b}_{i,j}&=& b_{i,j} - \alpha b_{i,k+2}b_{k,j},\;\; \quad j=1,\ldots,k+1, \quad i=k+1,\ldots,m, \label{eq:1stColproof:Coeffs2}\\
            \hat{c}_{j}  &=& c_{j}   - \alpha c_{k+2}  b_{k,j}, \;\; \qquad j=1,\ldots,k+1. \label{eq:1stColproof:Coeffs3}
        \end{eqnarray}
        \end{subequations}
        Then,
        \begin{equation}\label{eq:1stColproof:p_equal}
            \hat{p} = p,
        \end{equation}
        for any $k\in\left[1,\ldots,m\right]$ and $\alpha \in \CC$. \\      
\end{theorem}
\begin{proof}
Similarly to the proof of Theorem~\ref{thm:rowscaling}, we establish the following three statements:
\begin{subequations}\label{eq:1stColproof:Bcases}
\begin{eqnarray}
    \hat{Q}_{i+2}   &=& Q_{i+2} \quad \text{for } i = 1, \ldots, k-1, \label{eq:1stColproof:Bcases1}\\
    \hat{Q}_{k+2}   &=& Q_{k+2} + \alpha B_{k+2},\label{eq:1stColproof:Bcases2}\\
    \hat{Q}_i   &=& Q_i \quad \text{for } i = k+3, \ldots, m+2. \label{eq:1stColproof:Bcases3}
\end{eqnarray}
\end{subequations}
Together with the definition of $\hat{c}$ in \eqref{eq:1stColproof:Coeffs3},
this implies that the conclusion of the theorem in equation \eqref{eq:1stColproof:p_equal} holds.
Relation \eqref{eq:1stColproof:Bcases1} follows analogously to the proof in Theorem~\ref{thm:rowscaling}.

To prove \eqref{eq:1stColproof:Bcases2}, 
we insert the modified element $\hat{a}_{k,1}$ given in  \eqref{eq:1stcol_ABc} and observe that
\begin{eqnarray*}
    \hat{Q}_{k+2} &=& \left( (a_{k,1}+\alpha)Q_1 + a_{k,2}Q_2+\cdots+ a_{k,k+1}Q_{k+1}\right) B_{k+2} \\
                  &=& A_{k+2} B_{k+2} + \alpha B_{k+2} = Q_{k+2} + \alpha B_{k+2}.
\end{eqnarray*}

To prove \eqref{eq:1stColproof:Bcases3}, we prove that
$\hat{A}_{k+3}=A_{k+3},\ldots,\hat{A}_{m+2}=A_{m+2}$ and
$\hat{B}_{k+3}=B_{k+3},\ldots,\hat{B}_{m+2}=B_{m+2}$,
i.e., that the factors for $\hat{Q}_{k+3}$,\ldots, $\hat{Q}_{m+2}$ are unchanged. 
For $k+3$ we have 
\begin{equation} \label{eq:Afactor_proof} %
\begin{alignedat}{2}
      \hat{A}_{k+3} &= &&\hat{a}_{k+1,1} Q_1 + \cdots + \hat{a}_{k+1,k+1} Q_{k+1} + a_{k+1,k+2}\hat{Q}_{k+2} \\
              &= &&A_{k+3} - \alpha a_{k+1,k+2}( b_{k,1} Q_1 + \cdots + b_{k,k+1} Q_{k+1}) + a_{k+1,k+2} \alpha B_{k+2} \\
              &= &&A_{k+3} - \alpha a_{k+1,k+2} B_{k+2} + \alpha a_{k+1,k+2} B_{k+2} = A_{k+3},
\end{alignedat}
\end{equation}
where we have inserted \eqref{eq:1stColproof:Coeffs1} and \eqref{eq:1stColproof:Bcases2} to obtain the second equality.
The relation $\hat{B}_{k+3} = B_{k+3}$
can be shown analogously using \eqref{eq:1stColproof:Coeffs2} and \eqref{eq:1stColproof:Bcases2}. 
By induction, the corresponding relation for any factor.
Consequently, we have
\begin{equation}
    \hat{Q}_{i+2} = \hat{A}_{i+2} \hat{B}_{i+2} = A_{i+2} B_{i+2} = Q_{i+2}, \quad i = k+1,\ldots, m
\end{equation}
which proves \eqref{eq:1stColproof:Bcases3}. 

The conclusion \eqref{eq:1stColproof:p_equal}  follows from the same construction as in \eqref{eq:Afactor_proof}.
\end{proof}

\subsection{Normalized forms and unreduced schemes} 
The theorems in the previous section can be applied to impose a certain structure on the triplet $(A,B,c)$ without (or with very little) loss of generality. 
For any triplet $(A,B,c)$ we can invoke Theorem~\ref{thm:1stCol} repeatedly. By setting $\alpha$ to the negation of the element in the first column, we obtain matrices $A$ and $B$ with a first column containing only   zeros.  The first column corresponds to the addition of a scaled identity, which is independent (constant) with respect to the input $X$. Note that this can be done for any triplet, and no generality (in the sense of parameterized polynomials) 
is lost by assuming the first column is zero. 
\begin{definition}
    If $a_{1,1}=\cdots=a_{m,1}=b_{1,1}=\cdots=b_{m,1}=0$, we call the evaluation scheme a \emph{constant-free evaluation scheme}.
\end{definition}

\begin{corollary}
Any evaluation scheme is equivalent to a constant-free evaluation scheme.
\end{corollary}

Similarly, if we assume that the Hessenberg matrices $A$ and $B$ are unreduced \cite[p.~381]{Golub:2013:MATRIX}, i.e., the elements $a_{1,2},\ldots,a_{m,m+1}$ and $b_{1,2},\ldots,b_{m,m+1}$ are all nonzero, we can impose further structure by repeatedly applying Theorem~\ref{thm:rowscaling} with scaling determined by the corresponding last nonzero element of the row. This process imposes a normalization on each row.

\begin{definition}
A constant-free evaluation scheme satisfying $a_{1,2}=\cdots=a_{m,m+1}=b_{1,2}=\cdots=b_{m,m+1}=1$ is called a normalized evaluation scheme, and we call the triplet $(A,B,c)$ normalized.
\end{definition}

\begin{corollary}
Any evaluation scheme corresponding to $(A,B,c)$ where $A$ and $B$ are unreduced Hessenberg matrices is equivalent to a normalized evaluation scheme. 
\end{corollary}

\subsection{Further transformations}
For normalized constant-free evaluation schemes, we have $a_{11}=b_{11}=0$ and $a_{12} = b_{12} = 1$, meaning that the first multiplication always corresponds to squaring the input matrix, i.e.,
\begin{equation}\label{eq:Q2toQ3}
        Q_3 = X^2.
\end{equation}
This assumption can be made without loss of generality and will be used henceforth. This fact was already observed by Paterson and Stockmeyer \cite[p.~61]{Paterson:1973:OnTheNumber}.

We now derive further transformations under a technical assumption.
For the first two rows we assume that we have the structure of an unreduced constant-free triplet.\footnote{This assumption is made mostly to simplify the derivation, and similar results hold, e.g., when $a_{23}=0$.}
This in turn is equivalent to assuming $a_{23} = b_{23} = 1$. 

 Consider the following transformation  based on modifying the second row. We add $\alpha$ and $\beta$ to the second element in the second row in $A$ and $B$, respectively, i.e.,   
\begin{subequations}\label{eq:normalQ4hat}
\begin{eqnarray}
\hat{Q}_4 &=& ((a_{22}+\alpha)Q_2+Q_3)((b_{22}+\beta)Q_2+Q_3)\\    
 &=&Q_4+ (\alpha\beta + \alpha b_{22} + \beta a_{22} )X^2+(\alpha+\beta) X^3. \label{eq:normalQ4hatb}
\end{eqnarray}
\end{subequations}
If we enforce $\beta=-\alpha$, one term in the modification disappears, and we can compensate for the other term in a way such that  $\hat{Q}_5=Q_5,\ldots,\hat{Q}_m=Q_m$, i.e., they are unmodified. Consequently, the output polynomial is also unmodified.

\begin{theorem}\label{thm:b22zero}
Let $p$ be the polynomial associated with a constant-free triplet $(A,B,c)\in \CC^{m \times (m+1)} \times \CC^{m \times (m+1)} \times \CC^{m+2}$, satisfying $a_{23} = b_{23} = 1$, and let $\hat{p}$ be the polynomial associated with $(\hat{A},\hat{B},\hat{c})$ where
\begin{subequations}
    \begin{eqnarray}
\hat{A}&=&
\begin{bmatrix}
            0 & 1 &         &        &   &  &   \\
            0 & a+\alpha & a &        &     & &           \\
            0 & a & \hat{a}_{3,3} & a &     &&     \\
            0 & a & \hat{a}_{4,3} & a &  \ddots &      &     \\
            \vdots  & \vdots  &\vdots   & \vdots & \ddots&a &             \\
            0 & a & \hat{a}_{m,3} & a & \cdots & a   & a & 
\end{bmatrix}\\
\hat{B}&=&
\begin{bmatrix}
            0 & 1 &        &       &   &  &  \\
            0 & b-\alpha & b &       &     & &          \\
            0 & b & \hat{b}_{3,3} & b &     &&     \\
            0 & b & \hat{b}_{4,3} & b &  \ddots &      &     \\
            \vdots  & \vdots  &\vdots   & \vdots & \ddots&b & \\
            0 & b & \hat{b}_{m,3} & b & \cdots & b   & b & 
\end{bmatrix}
\end{eqnarray}
\end{subequations}
and
\begin{equation}
\hat{c} = 
\left[
    \begin{array}{c c c c c c}
        c_1 & c_2 & \hat{c}_3 & c_4 & \cdots & c_{m+2}
    \end{array}
\right].
\end{equation}
The modifications to the coefficients are given by
\begin{subequations}
    \begin{eqnarray}    
        \hat{a}_{i,3}&=&a_{i,3}-s a_{i,4} \label{eq:b22proof:coeffA}, \quad i = 3,\ldots,m. \\
        \hat{b}_{i,3}&=&b_{i,3}-s b_{i,4} \label{eq:b22proof:coeffB} \quad i = 3,\ldots,m.\\
        \hat{c}_3    &=&c_{3}  -s c_{4}, \label{eq:b22proof:coeffC}
    \end{eqnarray}
\end{subequations}
where
    \[ s = -\alpha^2 + \alpha (b_{22} - a_{22}) .\]
Then,
\begin{equation}
    p=\hat{p}. \label{eq:b22proof:peqp}
\end{equation}

\end{theorem}
\begin{proof}
    From the first row we have $\hat{Q}_3 = Q_3 = X^2$. The proof is based on establishing the following two statements:
    \begin{subequations}
        \begin{eqnarray}
        \hat{Q}_4 &=& Q_4 + s Q_3, \label{b22proof:Bcases1}\\
        \hat{Q}_{i+2} &=& Q_{i+2}, \quad i=3,\ldots,m. \label{b22proof:Bcases2}
        \end{eqnarray}
    \end{subequations}
    This, together with the definition of $\hat{c}$ in \eqref{eq:b22proof:coeffC} implies that the theorem conclusion \eqref{eq:b22proof:peqp} holds.

To prove \eqref{b22proof:Bcases1}, we use relation \eqref{eq:normalQ4hat}. More precisely, we substitute $\beta = -\alpha$ into \eqref{eq:normalQ4hatb} and obtain
\begin{equation}
    \begin{split}
        \hat{Q}_4 &= Q_4 + (-\alpha^2 + \alpha ( b_{22} - a_{22} ) ) X^2 = Q_4 + s Q_3.
    \end{split}
\end{equation}
To prove \eqref{b22proof:Bcases2}, we show that the 
factors for $Q_3,\ldots,Q_{m+2}$ are unchanged, i.e., 
$\hat{A}_3=A_3,\ldots, \hat{A}_{m+2}=A_{m+2}$
and 
$\hat{B}_3=B_3,\ldots, \hat{B}_{m+2}=B_{m+2}$.
For the first factor equality, we have
\begin{subequations}\label{eq:b22proof:A5}
\begin{eqnarray}
    \hat{A}_5 &=& a_{32}Q_2 + \hat{a}_{33}Q_3 + a_{34}\hat{Q}_4 \\
              &=& a_{32}Q_2 + a_{33} Q_3 - a_{34} s Q_3 + a_{34}Q_4 + a_{34} s Q_3 =         A_5,
\end{eqnarray}
\end{subequations}
where we have used \eqref{eq:b22proof:coeffA} and \eqref{b22proof:Bcases1} with $i=3$, in the second equality. The relation $\hat{B}_5 = B_5$ can be shown analogously using \eqref{eq:b22proof:coeffB} and \eqref{b22proof:Bcases1}. By induction, we can prove the corresponding factor relation for  $i=4,\ldots,m$. Consequently, we have
 $   \hat{Q}_{i+2} = \hat{A}_{i+2} \hat{B}_{i+2} = A_{i+2} B_{i+2} = Q_{i+2}$ for $i = 3,\ldots, m,$
which proves \eqref{b22proof:Bcases2}. The theorem conclusion \eqref{eq:b22proof:peqp} follows by the same construction as \eqref{eq:b22proof:A5}.
\end{proof}

For a given evaluation scheme satisfying $a_{23} = b_{23} = 1$, we can apply the previous theorem with $\alpha = b_{22}$, resulting in $\hat{b}_{22} = 0$. This shows that we can assume $b_{22} = 0$ for evaluation schemes that are unreduced in the first two rows, without loss of generality. 
For the second multiplication, under these assumptions, we get %
\begin{equation}\label{eq:b33proof:X4toX3}
    Q_4 = (a_{22}Q_2 + Q_3)Q_3 = a_{22}Q_2Q_3 + Q_3^2 = a_{22} X^3 + X^4.
\end{equation}

Next, we state and prove a theorem for the third row of the coefficient matrices. For this we assume that the first three rows are unreduced, i.e., $a_{23} = a_{34} = b_{23} = b_{34} = 1$.
The theorem is based on perturbing each of the free elements in the third row of $A$ and $B$, while simultaneously preserving the final output polynomial. In particular, we use perturbations with the following structure:
\begin{subequations}
    \begin{eqnarray*}
        \hat{a}_{32} &=& a_{32} + \alpha, \\
        \hat{b}_{32} &=& b_{32} - \alpha, \\
        \hat{a}_{33} &=& a_{33} + \beta, \\
        \hat{b}_{33} &=& b_{33} - \beta.
    \end{eqnarray*}
\end{subequations}
It turns out that when the perturbations have this particular structure, $Q_5$ is modified by the addition of a linear combination of $Q_3$, $Q_4$ and $X^3$. This is advantageous because we can compensate for $Q_3$ and $Q_4$, since we have access to these matrices directly. Moreover, we can ensure that the 
$X^3$-coefficient modification is zero by placing an additional
condition on  $\alpha$ and $\beta$, encoded in the equality $z(\alpha,\beta)=0$, where $z$ is a function explicitly given in the theorem.


\begin{theorem}\label{thm:b33adjust}
Let $p$ be the polynomial associated with the constant-free triplet $(A,B,c)\in \CC^{m \times (m+1)} \times \CC^{m \times (m+1)} \times \CC^{m+2}$ satisfying $a_{23} = a_{34} = b_{23} = b_{34} = 1$, and let $\hat{p}$ be the polynomial associated with $(\hat{A},\hat{B},\hat{c})$ where
\begin{subequations}
    \begin{eqnarray}
\hat{A}&=&
\begin{bmatrix}
            0 & 1 &        &       &   &  &  \\
            0 & a & a &       &     & &          \\
            0 & a+\alpha & a+\beta & a &    &&     \\
            0 & a & \hat{a}_{4,3} & \hat{a}_{4,4} &  \ddots &      &     \\
            \vdots  & \vdots  &\vdots   & \vdots & \ddots&a &            \\
            0 & a & \hat{a}_{m,3} & \hat{a}_{m,4} & \cdots & a   & a & 
\end{bmatrix}\\
\hat{B}&=&
\begin{bmatrix}
            0 & 1 &        &       &   &  &  \\
            0 & 0 & b &      &   & &          \\
            0 & b-\alpha & b-\beta & b &   & &     \\
            0 & b & \hat{b}_{4,3} & \hat{b}_{4,4} &  \ddots & & \\
            \vdots  & \vdots  &\vdots   & \vdots & \ddots&b & \\
            0 & b & \hat{b}_{m,3} & \hat{b}_{m,4} & \cdots & b   & b & 
\end{bmatrix}
\end{eqnarray}
\end{subequations}
and

\begin{equation}
\hat{c} = 
\left[
    \begin{array}{c c c c c c c}
        c_1 & c_2 & \hat{c}_3 & \hat{c}_4 & c_5 & \cdots & c_{m+2}
    \end{array}
\right].
\end{equation}
The modifications to the coefficients are given by
\begin{subequations}
    \begin{eqnarray}    
        \hat{a}_{i,3}&=&a_{i,3} - a_{i,5} s_1, \quad i = 4,\ldots,m, \label{eq:b33proof:coeffAk3}  \\
        \hat{a}_{i,4}&=&a_{i,4} - a_{i,5} s_2, \quad i = 4,\ldots,m, \label{eq:b33proof:coeffAk4}  \\
        \hat{b}_{i,3}&=&b_{i,3} - b_{i,5} s_1, \quad i = 4,\ldots,m, \label{eq:b33proof:coeffBk3}  \\
        \hat{b}_{i,4}&=&b_{i,4} - b_{i,5} s_2, \quad i = 4,\ldots,m, \label{eq:b33proof:coeffBk4}  \\
        \hat{c}_{3}&=&c_{3} - c_{5} s_1 \label{eq:b33proof:coeffck3}, \\
        \hat{c}_{4}&=&c_{4} - c_{5} s_2 \label{eq:b33proof:coeffck4},
    \end{eqnarray}
\end{subequations}
where
\begin{subequations}
    \begin{eqnarray}    
        s_1 &=&  \alpha ( b_{32} - a_{32} ) - \alpha^2,\label{b33proof:s1def}\\
        s_2 &=& \beta (b_{33} - a_{33}) - \beta^2. \label{b33proof:s2def}
    \end{eqnarray}
\end{subequations}
Then,
\begin{equation}\label{eq:b33proof:peqp}
    p=\hat{p}
\end{equation}
for any $\alpha$ and $\beta$ satisfying
\begin{equation}\label{eq:Zeq0}
    z(\alpha,\beta) \coloneqq \alpha(b_{33} - a_{33}) + \beta(b_{32} - a_{32}) - 2\alpha\beta - s_2a_{22} = 0.
\end{equation}
\end{theorem}



\begin{proof}
    It follows from the first two rows that $\hat{Q}_3 = Q_3 = X^2$, $\hat{Q}_4 = Q_4 = a_{22}X^3 + X^4$.
    The proof is based on establishing the following two statements:
    \begin{subequations}
        \begin{eqnarray}
        \hat{Q}_5 &=& Q_5 + s_1 Q_3 + s_2 Q_4 , \label{b33proof:Qcases1}\\
        \hat{Q}_{i+2} &=& Q_{i+2}, \quad i=4,\ldots,m. \label{b33proof:Qcases2}
        \end{eqnarray}
    \end{subequations}
    This, together with the definition of $\hat{c}$ in \eqref{eq:b33proof:coeffck3} and \eqref{eq:b33proof:coeffck4} implies that the theorem conclusion \eqref{eq:b33proof:peqp} holds.

    To prove \eqref{b33proof:Qcases1} 
we express $\hat{Q}_5$ in terms of the multiplication factors and substitute the definition of the modified coefficients. We obtain
    \begin{equation}\label{b33proof:Q5-1}
        \begin{split}
            \hat{Q}_5   &= \hat{A}_5\hat{B}_5 \\
                        &= ((a_{32} + \alpha)Q_2 + (a_{33} + \beta)Q_3 + Q_4)((b_{32} - \alpha)Q_2 + (b_{33} - \beta)Q_3 + Q_4) \\
                        &= \big(A_5 + (\alpha Q_2 + \beta Q_3) \big)\big(B_5 - (\alpha Q_2 + \beta Q_3) \big).    
        \end{split}
    \end{equation}
    This expression can be simplified by using $Q_5 = A_5 B_5$:
    \begin{equation}\label{b33proof:restform1}
    \begin{split}
        \hat{Q}_5   %
                    &= Q_5 + (B_5 - A_5)(\alpha Q_2 + \beta Q_3) - (\alpha Q_2 + \beta Q_3)^2.
    \end{split}
    \end{equation}
    The next step is to show that the difference between $Q_5$ and $\hat{Q}_5$ is a linear combination of $Q_3$, $Q_4$, and $X^3$, where the $X^3$-coefficient is given by $z(\alpha,\beta)$, which is zero by assumption. 
    We factorize the expression in \eqref{b33proof:restform1} to obtain
    \begin{equation}\label{b33proof:Q5-2}
    \begin{split}
        \hat{Q}_5   &= Q_5 + \big(B_5 - A_5 - (\alpha Q_2 + \beta Q_3)\big)\big(\alpha Q_2 + \beta Q_3\big).
    \end{split}
    \end{equation}
    Before proceeding we define
    \begin{equation}
        \begin{split}
            d_{32} := b_{32} - a_{32}, \\
            d_{33} := b_{33} - a_{33}.
        \end{split}
    \end{equation}
    This allows us to describe the difference between the multiplication factors more compactly:
    \begin{equation}
    \begin{split}
        B_5 - A_5 &= (b_{32} Q_2 + b_{33} Q_3 + Q_4) - (a_{32} Q_2 + a_{33} Q_3 + Q_4) \\
        &= (b_{32} - a_{32})Q_2 + (b_{33} - a_{33})Q_3.\\
        &= d_{32}Q_2 + d_{33}Q_3.
    \end{split}
    \end{equation}
    We substitute this expression into \eqref{b33proof:Q5-2} and simplify the first factor of the difference $\hat{Q}_5-Q_5$ such that 
    \begin{equation}
        \begin{split}
            \hat{Q}_5 &= Q_5 + \big( d_{32} Q_2 + d_{33}Q_3 -\alpha Q_2 -\beta Q_3 \big)\big(\alpha Q_2 + \beta Q_3\big) \\
            &=Q_5 + \big( (d_{32} - \alpha) Q_2 + (d_{33}-\beta)Q_3 \big)\big(\alpha Q_2 + \beta Q_3\big). \\
            \end{split}
            \end{equation}
    Next, we expand the final expression
        \begin{equation}
            \begin{split} %
                \hat{Q}_5   &= Q_5 + \alpha (d_{32} - \alpha) Q_2^2 + \big(\alpha(d_{33}-\beta) + \beta(d_{32} - \alpha)\big)Q_2Q_3 + \beta(d_{33}-\beta) Q_3^2 \\
                            &= Q_5 + s_1 Q_2^2 + (\alpha d_{33} + \beta d_{32} - 2\alpha\beta )Q_2Q_3 + s_2 Q_3^2, \\
            \end{split}
        \end{equation}
    where we have simplified in the last equality using \eqref{b33proof:s1def} and \eqref{b33proof:s2def}. By using $Q_3 = Q_2^2$, $Q_2Q_3 = X^3$ and $Q_3^2 = X^4$, we can rewrite this as
    \begin{equation}
        \hat{Q}_5 = Q_5 + s_1 Q_3 + (\alpha d_{33} + \beta d_{32} - 2\alpha\beta )X^3 + s_2 X^4.
    \end{equation}  
Finally, we use \eqref{eq:b33proof:X4toX3} in order to express $X^4$ in terms of $Q_4$ and $X^3$
    \begin{equation}
        \begin{split}
            \hat{Q}_5 &= Q_5 + s_1 Q_3 + (\alpha d_{33} + \beta d_{32} - 2\alpha\beta )X^3 + s_2 (Q_4 - a_{22}X^3) \\
            &= Q_5 + s_1 Q_3 + \big( \alpha d_{33} + \beta d_{32} - 2\alpha\beta - s_2 a_{22} \big)X^3 + s_2 Q_4 \\
            &= Q_5 + s_1 Q_3 + z(\alpha,\beta)X^3 + s_2 Q_4 \\
            &= Q_5 + s_1 Q_3 + s_2 Q_4,
        \end{split}
    \end{equation}
    where we have used that $z(\alpha,\beta) = 0$ in the last equality. This proves statement \eqref{b33proof:Qcases1}.

    To prove \eqref{b33proof:Qcases2} we show that the multiplication factors are unchanged, i.e., $\hat{A}_{6}=A_{6},\ldots\hat{A}_{m+2}=A_{m+2}$ and $\hat{B}_{6} = B_{6},\ldots \hat{B}_{m+2} = B_{m+2}$. For the first factor we have
\begin{equation} \label{eq:b33proof:A6}
    \begin{split}
        \hat{A}_{6} &= (a_{42}Q_2 + \hat{a}_{43}Q_3 + \hat{a}_{44}Q_4 + a_{45}\hat{Q}_5) \\
                    &=  a_{42} Q_2 + (a_{43} - a_{45}s_1) Q_3 + (a_{44} - a_{45}s_2) Q_4 + a_{45} Q_5 + a_{45}( s_1 Q_3 + s_2 Q_4 ) \\
                    &= A_6 - a_{45} s_1 Q_2 - a_{45} s_2 Q_4 + a_{45}( s_1 Q_3 + s_2 Q_4 ) \\
                    &= A_6,
    \end{split}
\end{equation}
where we have used \eqref{eq:b33proof:coeffAk3}, \eqref{eq:b33proof:coeffAk4} and \eqref{b33proof:Qcases1} in the second equality. 
The relation $\hat{B}_6 = B_{6}$ follows analogously using \eqref{eq:b33proof:coeffBk3}, \eqref{eq:b33proof:coeffBk4} and \eqref{b33proof:Qcases1}.

The corresponding relation can be shown for all factors using induction. Consequently, we have
\begin{equation}
    \hat{Q}_{i+2} = \hat{A}_{i+2} \hat{B}_{i+2} = A_{i+2} B_{i+2} = Q_{i+2},\quad i = 4,\ldots, m,
\end{equation}
which proves \eqref{b33proof:Qcases2}. The theorem conclusion \eqref{eq:b33proof:peqp} follows by the same construction as \eqref{eq:b33proof:A6}.
\end{proof}


\subsection*{Free variables in Theorem~\ref{thm:b33adjust}}
Note that Theorem~\ref{thm:b33adjust} includes a scalar-valued condition, $z(\alpha,\beta) = 0$, involving two scalar variables.
Let $d_{32}=b_{32}-a_{32}$ and $d_{33}=b_{33}-a_{33}$ 
be defined as in the proof of the theorem.
If we let $\alpha$ be given, we get a quadratic equation in $\beta$:
\begin{equation}
    d_{33}\alpha  = a_{22}\beta^2 + \beta\left(2\alpha + d_{32} - a_{22}d_{33}\right).
\end{equation}
When $a_{22}\neq 0$,  
the solution $\beta$ to the equation is explicitly available from the solution of the quadratic equation. This has a disadvantage of introducing a square root operation.
In a real setting, this can yield complex coefficients.

This is related to the result in paper \cite{Sastre:2018:SMARTPOLY} as follows. The results suggest several approaches to evaluate
polynomials with a low number multiplications. For the case $m=3$, the formulas
 \cite[Eq. (31)]{Sastre:2018:SMARTPOLY} involve a square root, 
 and indeed that approach can be derived from the above transformation 
 with $\alpha=-a_{32}$ and solving for $\beta$.

Suppose $\beta$ is given. Then, the solution for $\alpha$ can be expressed as
\begin{equation}
    \alpha = \frac{a_{22}\beta^2 + (d_{32} - a_{22}d_{33})\beta}{
    d_{33} - 2\beta}
\end{equation}
with the condition $\beta \neq \frac{1}{2}d_{33}$. To reframe this condition in
terms of entries in matrices $A$ and $B$, we use a change of variables
that essentially generalizes \cite{Rubebsson:2025:SP8} (where it is given for $m=3$):
\begin{equation}
    \beta = \frac{1}{2}d_{33} + r, \quad r \neq 0.
\end{equation}
With this choice, we obtain updated table entries:
\begin{equation}
\begin{split}
    \hat{a}_{33} &= a_{33} + \beta = a_{33} + \frac{b_{33}-a_{33}}{2} + r = \frac{a_{33}+b_{33}}{2} + r, \\
    \hat{b}_{33} &= b_{33} - \beta = b_{33} - \frac{b_{33}-a_{33}}{2} - r = \frac{a_{33}+b_{33}}{2} - r.
\end{split}
\end{equation}
This leads to the relation
\begin{equation}
    \hat{a}_{33} - \hat{b}_{33} = 2r.
\end{equation}
Consequently, any evaluation scheme that is
unreduced in the first three rows can be transformed into an equivalent scheme satisfying
\begin{equation}
    b_{33} = a_{33} + 1,
\end{equation}
since we can choose any nonzero $r$. Thus, when considering evaluation schemes that are unreduced in the first three rows, we can assume without loss of generality that $b_{33} = a_{33} + 1$. This assumption allows us to reduce the number of variables without introducing additional constraints, as explained in the following section.

\section{Minimality of parameterization}\label{sec:minimality}

\subsection{Minimality of unreduced evaluation schemes}
As a consequence of the equivalence theorems in the previous
section, we can conclude that any unreduced evaluation scheme can be assumed to have the form

\begin{equation} \label{eq:parameterization}
[A|B]=\left[
\begin{array}{ccccccc|cccccccc}
0 & 1 &   &   &   &   && 0 & 1 &   &   &   &   \\
0 & a & 1 &   &   &   && 0 & 0 & 1 &   &   &   \\
0 & a & a_{33} & 1 &   &   && 0 & b & {a_{33}+1} & 1 &   &   \\
0 & a & a & a & 1 &   && 0 & b & b & b & 1 &   \\
\vdots & \vdots & \vdots & \vdots & \ddots & \ddots && 
\vdots & \vdots & \vdots & \vdots & \ddots & \ddots\\
0 & a & a & a & \cdots & a & 1 & 0 & b & b & b & \cdots & b & 1 \\
\end{array}\right],
\end{equation}
with an arbitrary $c$-vector. The dimension of $\pistar{m}$ 
is determined from the rank of the Jacobian for a generic element in $\pistar{m}$. A generic element in $\pistar{m}$ is unreduced; hence, the dimension is bounded by the number of free variables in a parameterization of the unreduced schemes.
The parameterization \eqref{eq:parameterization} contains $m^2$ parameters and, therefore, similar to \eqref{eq:dim_s_bound}, we have the following.
\begin{corollary}\label{thm:dim_2m}
For $m>2$, we have
\[\dim(\Pi^*_{2^m})\le m^2.\]
\end{corollary}

It turns out that this bound is sharp; i.e., we can establish a lower bound which is also $m^2$. To this end it is sufficient to find a triplet $(A,B,c)$ such that the Jacobian of the output $p(X)$ with respect to the free variables
has rank $m^2$. 
Many trivial choices of unreduced $(A,B,c)$, e.g., a triplet corresponding to only a sequence of squaring operations, lead to a points in 
$\pistar{m}$ whose rank is less than $m^2$. 
The triplet with a full Jacobian rank and leading to the simplest formulas for the Jacobian that we could find is the following: 
\begin{equation}\label{eq:proof_example}
     [A|B]=\left[
\begin{array}{cccccc|cccccc}
0 & 1 &   &   &       &&  0 & 1 &   &   &   &   \\
0 & 1 & 1 &   &       &&  0 & 1 & 1 &   &   &   \\
0 & 1 & 0 & 1 &       &&  0 & 1 & 1 & 1 &   &   \\
\vdots & \vdots & \vdots & \ddots & \ddots&&  \vdots & \vdots & \vdots & \ddots & \ddots &   \\
0 & 1 & 0 &  \cdots & 0 & 1  &  0 & 1 & 1 &  \cdots & 1 & 1 \\
\end{array}
\right] 
\end{equation}
and $c=[0,\ldots,0,1]$. 

In order to bound the rank, it is sufficient to find $m^2$ linearly independent columns in the Jacobian. 
A sufficient condition for linear independence for polynomials is that they have distinct degrees. 
Based on that reasoning, we now establish linear combinations of partial derivatives of the output of the evaluation, i.e., $p(X)$,  with respect to the elements of $A,B$ and $c$. We obtain $m^2$ distinct degrees.

Although the choice \eqref{eq:proof_example} leads to the simplest derivation we could establish analytically, it is admittedly not very simple. The proof of the general case can be found in the appendix. For illustration, we sketch the derivation for $m=4$, specifying the partial derivatives, which can be computed with symbolic computation tool (e.g., those described in the next section). In the full proof in the appendix we provide the details without such tools. We base the proof, as well as this sketch on a separation into cases, each case leading to a polynomial degrees that are distinct, and adding up to $m^2=16$ different degrees. 

Case 1 corresponds to forming derivatives with respect to elements of $c$:
$\frac{\partial p}{\partial c_1} = 1,\ 
\frac{\partial p}{\partial c_2} = X,\ 
\frac{\partial p}{\partial c_3} = X^2,\ 
\frac{\partial p}{\partial c_4} = Q_4,\ 
\frac{\partial p}{\partial c_5} = Q_5,\ 
\frac{\partial p}{\partial c_6} = Q_6$.
Since $(A,B,c)$ in \eqref{eq:proof_example} is unreduced, $Q_4,Q_5$ and $Q_6$ have maximal degree and we have that 
\begin{equation} \label{eq:deg_ci}
\deg\left(\frac{\partial p}{\partial c_i}\right)=
\deg(Q_i)=2^{i-2}, i=2,\ldots,m+2.
\end{equation}
Case 2 corresponds to forming derivatives with respect to $a_{i,j}$. For example, we have

\begin{eqnarray}
\frac{\partial p}{\partial a_{4,4}}&=&X^{12} + 6X^{11} + 16X^{10} + 26X^{9} + \mathcal{O}(X^{8}) \label{eq:case2_a44} \\
\frac{\partial p}{\partial a_{4,2}}&=&X^{9} + 4X^{8} + 7X^{7} + 8X^{6} + 7X^5 + \mathcal{O}(X^{4}).\label{eq:case2_a42}
\end{eqnarray}

As shown in the proof of the theorem, the general formula for the degrees in Case~2 is
\begin{equation}\label{eq:case2_degformula}
     \deg\left(\frac{\partial p}{\partial a_{i,j}}\right)=2^{m}-2^{i-1}+2^{j-2}, \textrm{ for }i=2,\ldots,m,\;\;j=2,\ldots,i
\end{equation}

Case 3 stems from the observation that we obtain the same degree if we differentiate with respect to $a_{i,j}$ and $b_{i,j}$. 
Therefore we form the difference in order to obtain a different degree. 
For example, 
\begin{equation}\label{eq:case3_b42}
\frac{\partial p}{\partial b_{4,2}}=X^{9} + 4X^{8} + 7X^{7} + 8X^{6} + 6X^5+\mathcal{O}(X^{4}) 
\end{equation}
and  the difference with \eqref{eq:case2_a42} is
\begin{equation}\label{eq:case3_x42}
\frac{\partial p}{\partial b_{4,2}}-\frac{\partial p}{\partial a_{4,2}}=
 -X^{5}  -2X^{4}  -2X^{3}.
\end{equation}

The degree $5$ is distinct from the degrees in Cases 1 and 2. In the appendix we prove that the general formula for the degrees in Case 3 is
\begin{equation}\label{eq:case3_degformula}
    \deg\left(\frac{\partial p}{\partial b_{i,j}}-
\frac{\partial p}{\partial a_{i,j}}\right) = 2^m-2^i+2^{i-2}+2^{j-2},
\textrm{ for }i=3,\ldots,m,\;\;j=2,\ldots,i-1
\end{equation}
and that they are distinct from previous cases.

\newpage

Case 4 is based on the observation that if we use the idea from Case 3 for $j=i=m$, we get a degree which already included in Case 1. In order to establish a degree not present in the previous cases, we must form a linear combination of partial derivatives with respect to the elements $c_{m+1}$, $a_{m,m}$, $b_{m,m}$, $a_{m,2}$ and $b_{m,2}$.
Consider the following identities:
\begin{eqnarray}
\frac{\partial p}{\partial b_{4,4}}&=&X^{12} + 6X^{11} + 16X^{10} + 26X^{9} + \mathcal{O}(X^{8}) \\
\frac{\partial p}{\partial b_{4,4}}-\frac{\partial p}{\partial a_{4,4}}&=&-X^{8} - 4X^{7} - 7X^{6} - 6X^{5} + \mathcal{O}(X^{4}) \label{eq:case4_x44}\\
\frac{\partial p}{\partial c_{5}}&=&X^{8} + 4X^{7} + 7X^{6} + 8X^{5} + \mathcal{O}(X^{4}).\label{eq:case4_c5}
\end{eqnarray}
By forming the sum of equations \eqref{eq:case4_x44} and \eqref{eq:case4_c5}, we obtain 
\[\frac{\partial p}{\partial b_{4,4}}-\frac{\partial p}{\partial a_{4,4}}+\frac{\partial p}{\partial c_{5}}= 2X^{5} + 4X^{4} + 3X^{3} + X^{2}. \]
Noting from \eqref{eq:case3_x42} that this degree coincides with Case 3 for $i=4$ and $j=2$,
and we can reduce the degree by forming the sum 
\[
    \frac{\partial p}{\partial b_{4,4}}
    -\frac{\partial p}{\partial a_{4,4}}
+
    \frac{\partial p}{\partial c_{5}}
+
    2\left(
        \frac{\partial p}{\partial b_{4,2}}-
        \frac{\partial p}{\partial a_{4,2}}
    \right)
=
-X^3+X^2.
\]
The resulting degree 3 is distinct from those in Cases 1-3.  

From the reasoning above, in this example, we have 6 distinct degrees from Case~1, 6 distinct degrees from Case~2, 
3 distinct degrees from Case~3 and 1 degree from Case~4, which yield a total of $m^2=16$ distinct degrees.

In the general case, i.e., when $m>4$, we also provide a Case 5, leading to an additional to $m-4$ distinct degrees. The idea for this case is very similar to that of Case 4, and corresponds to forming a linear combination of partial derivatives with respect to the elements $a_{i,i}$, $b_{i,i}$ $a_{i+1}$, $a_{i,2}$ and $b_{i,2}$ for $i<m$.

\begin{theorem}\label{thm:dim_2m_conj}
For $m>2$, we have
\begin{equation}\label{eq:dim_2m_eq}
\dim(\Pi^*_{2^m})= m^2.
\end{equation}
\end{theorem}
\begin{proof}
    See Section~\ref{sec:appendix_proof}.
\end{proof}

As a consequence of Theorem~\ref{thm:dim_2m_conj}, the parameterization \eqref{eq:parameterization} is minimal.
To our knowledge, this is the first minimal parameterization of unreduced schemes in $\pistar{m}$.

\subsection{Reduced evaluation schemes}
In practice, unreduced evaluation schemes are rarely useful for evaluating a given polynomial of high degree $d$, because of the limited number of degrees of freedom in $\pistar{m}$ for large $m$, in comparison to  $\dim(\pideg{2^m})=d+1=2^m+1$. If the Hessenberg matrices $A$ and $B$ are reduced, we obtain output polynomials of lower degree. For example the pair
\begin{equation} \label{eq:reduced_sixmult}
[A|B]=\left[
\begin{array}{ccccccc|cccccccc}
0 & 1 &   &   &   &   && 0 & 1 &   &   &   &   \\
0 & \times & 1 &   &   &   && 0 & 0 & 1 &   &   &   \\
0 & \times & \times & 1 &   &   && 0 & \times & \times & 1 &   &   \\
0 & \times & \times & \times & 1 &   && 0 & \times & \times & 1 & 0 &   \\
0 & \times & \times & \times & \times & 1 && 0 & \times & \times & \times & 1 & 0 \\
0 & \times & \times & \times & \times & \times & 1 & 0 & \times & \times & \times & \times & 1 & 0 \\
\end{array}\right]
\end{equation}
corresponds to $m=6$ multiplications but results in a polynomial of degree $32$.
By a single reduction, we mean setting the last nonzero element in a row to zero in either $A$ or $B$. With $r$ reductions we mean the repeated application of a single reduction. Note that we can still normalize each row since the transformation theorems are also applicable to reduced systems. Hence, we lose one degree of freedom with every reduction, and due to Corollary~\ref{thm:dim_2m}, we expect the corresponding dimension to be
\begin{equation} \label{eq:m2_r_limit}
    m^2-r.
\end{equation}
In the following section we proceed by studying reduced evaluation schemes.
More precisely, we study reduced evaluation schemes 
that lead to specific polynomial degrees and describe ways to compute $(A,B,c)$ for that reduction structure for a polynomial given in a monomial basis.

\begin{figure}[h]
  \begin{center}
    \subfigure[]{\scalebox{0.9}{
\input{gfx/reduced1.tex}  
    }}%
    \subfigure[]{\scalebox{0.9}{
\input{gfx/reduced2.tex}  
    }}
    \caption{The polynomial degrees when the Hessenberg matrices in the triplet $(A,B,c)$ are reduced Hessenberg matrices. 
    The parameter $\epsilon=0.1$ is selected for visualization purposes. Any point above the blue curve can be discarded to not completely contain the corresponding polynomial subset, in the sense of \eqref{eq:maxproblem}, due to an insufficient number of degrees of freedom.
      \label{fig:reduced}
    }
  \end{center}
\end{figure}

\section{Polynomial subsets}\label{sec:casestudies}

This section is devoted to the study of the question: \emph{What is the 
    largest $d$ such that all $d$-degree polynomials can be computed with $m$ multiplications?}
Formally, we use two versions of this problem
    \begin{equation}\label{eq:maxproblem}    
    \max\{d: \Pi_d \subset \Pi_{2^m}^*\}
    \leq     \max\{d: \Pi_d \subset \widebar{\Pi_{2^m}^*}\}.
    \end{equation}    
    We consider the left-hand side when possible, and otherwise study the right-hand side in order to avoid limit cases similar to \eqref{eq:epsilonexample}.
     
The previous section stressed the use of reduced matrices. For a given reduction, we can compute the corresponding polynomial degree. 
Hence, we can investigate candidate solutions to \eqref{eq:maxproblem} by considering all possible reductions. 
This approach is depicted in Figure~\ref{fig:reduced}, which illustrates all combinations of reductions for up to $m \leq 7$ multiplications and $r \leq 7$ reductions.
For instance, \eqref{eq:reduced_sixmult} corresponds to the scenario (horizontal axis) with $m=6$ multiplications and $r=3$ reductions and achieves a polynomial degree of $32$, as shown on the vertical axis of the figure.

To identify optimal solutions to \eqref{eq:maxproblem}, higher polynomial degrees are advantageous.
However, excessively high degrees may result in insufficient degrees of freedom.
More precisely,  \eqref{eq:m2_r_limit}  gives a bound on the degree of the 
admissible candidate solutions to \eqref{eq:maxproblem} for a given number of multiplications and reductions.
This is visualized with a  blue line in the figure.
Hence, for the purpose of studying \eqref{eq:maxproblem} we can disregard points above this line.

The problem \eqref{eq:maxproblem} becomes increasingly complex as $m$ increases.  For  $m=3$, the problem is essentially already solved in \cite{Sastre:2018:SMARTPOLY} since an explicit procedure is provided to compute almost all polynomials of degree $8$ with $m=3$ multiplications, i.e., 
\[
 \max\{d: \Pi_d \subset \widebar{\Pi_{2^3}^*}\}=8.
\]
For $m=4$, similar constructions are also provided in \cite{Sastre:2018:SMARTPOLY}, yielding a method for degree $12$.
An alternate approach requiring fewer assumptions on the monomial coefficients is given in Section~\ref{sec:fourmult}.
From Figure~\ref{fig:reduced}, we see that this is the highest admissible degree and therefore conclude that it is optimal.

For $m\ge 5$ we were not able to solve the problem analytically and instead resorted to computational tools. More precisely, we frame the problem with a given reduction and structure as a system of polynomial equations. For $m=5$, we use the package HomotopyContinuation.jl \cite{HomotopyContinuation.jl} to find, as far as we can tell,  all solutions. For $m=6$ and $m=7$, the system was too complicated, and we were not able to find solutions with this package. However, by using the software \cite{jarlebring:2023:graphs}, we were able to construct locally convergent iterative methods and successfully found solutions in all specified test cases. These simulations combined with reasoning based on admissible degrees in Figure~\ref{fig:reduced}, led us to state conjectures concerning the solution to \eqref{eq:maxproblem}. For reproducibility, all simulations (including starting values) are provided in the publicly available GitHub repository: \url{https://github.com/GustafLorentzon/polynomial-set-paper}.


\subsection{Four multiplications}\label{sec:fourmult}
When we study $m=4$, we identify from Figure~\ref{fig:reduced} that the highest degree of the admissible polynomials is $12$, since without reductions we only have 
$16$ degrees of freedom, which cannot parameterize $\pideg{16}$.  The solution to \eqref{eq:maxproblem} is indeed 12. This can already be concluded from the method in \cite[p.~237]{Sastre:2018:SMARTPOLY} which is a method to evaluate polynomials of degree 12, given in their monomial basis, using only $m=4$ multiplications.
In our terminology this corresponds to the reduction $a_{2,3}=0$ and, additionally, $a_{4,2}=0$.  The method in
\cite[p.~237]{Sastre:2018:SMARTPOLY} involves 
square roots of expression containing the monomial coefficients, 
roots of a polynomial of degree four,
as well as divisions of certain quantities leading to exceptions corresponding to some limit cases. Therefore, from \cite{Sastre:2018:SMARTPOLY} we conclude that 
\begin{equation}
\max\{d: \Pi_d \subset \pistarbar{4}\}=12
\end{equation}
in a complex sense.

We now present a slightly more general method for evaluating any polynomial in $\pideg{12}$ with four multiplications. The new method does not involve square roots or divisions except for the leading monomial coefficient. Consider the evaluation schemes with the following structure
\begin{subequations}\label{eq:ABfourmult}    
\begin{align}
        \left[\begin{array}{c | c}
                A & B
            \end{array}\right]
            &= \left[\begin{array}{c c c c c | c c c c c}
            0 & 1 &     &           &           & 0 & 1      &                 & \\
            0 & 1 & 0 &        &           & 0 & 0      & 1               & \\
            0 & a_{32} & a_{33} & 1 &           & 0 & 0& 0& 1\\
            0 & a_{42} & a_{43} & a_{44} & 1    & 0 & b_{42} & b_{43} & a_{44}+1 & 1
            \end{array}\right] \\
            c &=
            \left[
            \begin{array}{cccccc} 
            c_1 & c_2 & c_3 & c_4 & c_5
            \end{array}
            \right].
\end{align}
\end{subequations}
Suppose $\alpha_0,\ldots,\alpha_{12}$ represent a given polynomial $p(X) = \alpha_0 I + \alpha_1 X + \cdots + \alpha_{12} X^{12} \in \pideg{12}$. If we expand the parameterization, we obtain a multivariate polynomial system of equations---one for each monomial coefficient in the output polynomial. In the terminology of Section~\ref{sec:variety}, the system corresponds to considering the $0$th,$\ldots$,$12$th derivatives of the equation $\Phi(A,B,c)(x)=p(x)$ with respect to $x$, evaluated at $x=0$.
In this case, we have $13$ equations in $13$ variables. To solve solve this system,
we first introduce the auxiliary variables 
$\beta_{43}= b_{43} + a_{43}$ and $\beta_{42}=b_{42} + a_{42}$.
This system is explicitly solvable by considering the equations in the output polynomial ordered 
in descending degree, so that the equation corresponding to $\alpha_{12}$ is treated first.  In this sense, the system is triangular, which was also crucial for the construction in \cite[Section 3]{Sastre:2018:SMARTPOLY}.
The solution to the system is given by the following sequence of equations:
\begin{subequations}
    \begin{align}
            c_6 &= \alpha_{12} \\
            a_{33} &= \frac{1}{2}\left(\frac{\alpha_{11}}{c_6}\right) \\
            a_{32} &= \frac{1}{2}\left( \frac{\alpha_{10}}{c_6} - a_{33}^2\right) \\
            a_{44} &= \frac{1}{2}\left( \frac{\alpha_{9}}{c_6} - 2a_{32}a_{33} - 1\right) \\
            \beta_{43} &= \frac{\alpha_{8}}{c_6} - \left( a_{33} + 2a_{33}a_{44} + a_{32}^2 \right) \\
            \beta_{42} &= \frac{\alpha_7}{c_6} - \left( a_{32} + a_{33}\beta_{43} + 2a_{32}a_{44} \right) \\
            c_5 &= \alpha_6 - c_6\left(a_{44} + a_{44}^2 + a_{33}\beta_{42} + a_{32}\beta_{43} \right) \\
            a_{43} &= \frac{c_5}{c_6} - \left(a_{33}\frac{c_5}{c_6} + a_{44}\beta_{43} + a_{32}\beta_{42}\right) \\
            a_{42} &= \frac{\alpha_4}{c_6} - \left(a_{32}\frac{c_5}{c_6} + a_{44}\beta_{42} + a_{43}\beta_{43} - a_{43}^2\right) \\
            c_4 &= \alpha_3    - c_6\left(a_{43}\beta_{42} + a_{42}\beta_{43} - 2a_{42}a_{43}\right) \\
            c_3  &= \alpha_2    - c_6\left(a_{42}\beta_{42} - a_{42}^2\right).
            \end{align}
\end{subequations}
With the conditions $c_2 = \alpha_1$, 
$c_1  = \alpha_0$,
$b_{43} = \beta_{43} - a_{43}$ and 
$b_{42} = \beta_{42} - a_{42}$, we have explicitly computed all variables in \eqref{eq:ABfourmult}.
    
Recall that $\alpha_{12}\neq 0$ for $p\in \pideg{12}$; therefore, we have made no assumptions other than the degree of the polynomial. Moreover, the formulas preserve the algebraic structure of the variables, e.g., if $\alpha_0,\ldots,\alpha_{12}\in \RR$, then $(A,B,c)$ is a real triplet, so the evaluation  coefficients are real.  We conclude that
\begin{equation}
\max\{d: \Pi_d \subset \pistar{4}\}=12
\end{equation}
holds in both a real and a complex sense.

\subsection{Five multiplications}\label{sec:fivemult}
For $m=5$ multiplications we see in Figure~\ref{fig:reduced} that the highest
degree of the admissible polynomials is $20$. 
To our knowledge, the state of the art is $d=18$ as given in \cite[Equation (17)-(19)]{saib21} with $s=2$ which is based on  \cite{Sastre:2018:SMARTPOLY} combined
with the Paterson--Stockmeyer evaluation. In our terminology, that approach corresponds to the reduction $a_{23}=a_{56}=0$ and additionally imposing $a_{42}=0$. The reduction in Figure~\ref{fig:reduced} leading to a polynomial of degree $d=20$ corresponds to $a_{45}=a_{56}=0$. 

We were not able to explicitly derive a solution to the multivariate polynomial system for the structure with $d=20$ analytically; instead, we needed to resort to computational tools. The Julia package HomotopyContinuation.jl \cite{HomotopyContinuation.jl} includes methods to solve polynomial systems of equations based on numerical continuation and with advanced initialization of starting points for the homotopy method. 
We implemented the system equations derived from considering each monomial coefficient for the data structures of this package. Since we have more variables than equations for this structure, we empirically fixed some variables, choosing which variable to fix by trying to reduce the total degree of the system as much as possible. We additionally solved those equations that could be solved explicitly, e.g., the first and last equations.

With this setup, we were able to compute $(A,B,c)$ for a large number of given polynomials $p\in\pideg{20}$, including the truncated Taylor expansion of the exponential as well as the function $1/(1-x)$. 
The simulations were done in Julia, and the code is given in the GitHub repository. Moreover, code to actually evaluate polynomials is provided in both Julia and Matlab. 
For conciseness, we report only the numbers for the matrix exponential in the following. 

In an attempt to prevent large condition numbers, we sought solutions that did not have excessively large or small values. In this case, we mitigated large numbers by scaling the input.  For the exponential, we approximate $e^{\alpha x}$ with $\alpha=8$ fixed, since this makes $c_7=\alpha^{20} / 20!\approx 0.47$ in the order of magnitude one. Although the input scaling can be reversed by transforming entries in the table, it was not deemed numerically useful and therefore it was not included in the following presentation of results.

HomotopyContinuation.jl found several solutions and the solution with the smallest values was the following
\begin{subequations}\label{eq:five_mult_AB}
\begin{align}
[A|B]&=\left[
\begin{array}{cccccc|cccccc}
0 & 1 &   &    &    &    & 0 & 1   &    &    &    &    \\
0 & 0 & 1 &    &    &    & 0 & \frac12 & 1  &    &    &    \\
0 & 0 & 2 & 1  &    &    & 0 & b & 1 & 1  &    &    \\
0 & a & a & 1  & 0  &    & 0 & b & b & b & 1  &    \\
0 & a & a & a  & 1  & 0  & 0 & b & b & b & b & 1  
\end{array}
\right] \\
c &=\left[
            \begin{array}{cccccccc} 
            c & c & c & c & c & c & c
            \end{array}
            \right]
\end{align}
\end{subequations}
where the missing values are given in Table~\ref{fig:five_missing}. The Jacobian of the polynomial system evaluated at this solution has a condition number of $8.1\cdot10^{2}$.
\begin{table}[h]
  \begin{center}
    \begin{tabular}{|c|c|c|c|}

$a_{42}$ & 2.3374451754385963 & $c_{1}$ & $1$\\ 	 
$a_{43}$ & $-41/16 $ 	       & $c_{2}$ & $\alpha$\\
$a_{52}$ & 2.8309861554443847 & $c_{3}$ & -6.657689892163032 \\
$a_{53}$ & 8.7616118485412    & $c_{4}$ & 50.445902670306005 \\		
$a_{54}$ & 5.123957592622475  & $c_{5}$ & 19.754729172913187 \\	
$b_{32}$ & 1.4484649122807018 & $c_{6}$ & 2.8090057997411706 \\
$b_{42}$ & 6.389966463262669  & $c_{7}$ & $\alpha^{20}/20!$ \\ 
$b_{43}$ & 6.697361614351532 \\
$b_{44}$ & 2.1451472591988834 \\
$b_{52}$ & -2.458444697550387 \\
$b_{53}$ & -3.6724346694235876 \\
$b_{54}$ & 16.044090747953085 \\
$b_{55}$ & 7.557067023178642 \\

    \end{tabular}
    \caption{Non-specified values in \eqref{eq:five_mult_AB}
      \label{fig:five_missing}
    }
  \end{center}
\end{table}

Since we were able to find solutions in the case studies we conjecture that this 
corresponds to a realization of a method for the maximum polynomial degree.
\begin{conjecture}\label{conj:fivemult}
  \[
\max\{d: \Pi_d \subset \pistarbar{5}\}=20.
  \]
\end{conjecture}
The upcoming work \cite{sastre2025advancing} suggests that there is indeed a constructive way to form such evaluation schemes for five multiplications.

\subsection{Six multiplications}
To our knowledge, the state of the art for $m=6$ multiplications
is $d=24$,  again given in 
\cite{saib21} with $s=3$. In our terminology, that corresponds to
the reduction $a_{23}=a_{33}=a_{34}=a_{67}=0$. Similar to the situation
for $m=5$, this is not  the highest admissible degree. From Figure~\ref{fig:reduced} we see that the highest admissible degree is $d=32$.
With $r=3$, we can use the reduction \eqref{eq:reduced_sixmult}.

Unfortunately, the application of the package HomotopyContinuation.jl was not successful for this case. We have $33$ equations, and $m^2-r=33$ unknowns.
The creation of initial vectors for the homotopy continuation seems too 
computationally demanding, 
likely related to the high total degree 
of the polynomial system. Instead we used the package GraphMatFun.jl \cite{jarlebring:2023:graphs} to create a locally convergent iterative solution method. The system is highly ill-conditioned and a standard Newton approach was not successful.
Instead, we employ an iteratively regularized Newton's method. Following the approach in \cite{eriksson1996regularization}, we compute a \emph{Tikhonov--Newton step} by applying Newton's method to a Tikhonov-regularized system.
The Tikhonov--Newton step can be computed using the singular value decomposition of the Jacobian matrix, which is directly available from the graph representation in GraphMatFun.jl. The best results were obtained by using Armijo step-length damping and selecting between a Newton and a Tikhonov--Newton step in a greedy fashion.
For the sake of reproducibility, the starting vectors are given explicitly in the software available in the GitHub repository. In order to determine conclusively that a solution is found, the solution was post-processed with high-precision arithmetic (BigFloat in Julia) so that the first 50 decimals of the coefficients appear accurate.  This was done for this particular simulation as well as for all subsequent simulations. 

The above approach with the structure \eqref{eq:reduced_sixmult} was successful in finding a solution vector for the problem corresponding to the Taylor expansion of the matrix exponential, \emph{using complex arithmetic}. Unfortunately, we were not able to find a real coefficient vector. 

From an application viewpoint, real coefficients are advantageous, e.g., since the evaluation of $p(A)$ can be done in real arithmetic if $A$ is real.
In order to find a real evaluation scheme, we investigated instead the  
structure
\[
[A|B]=\left[
\begin{array}{ccccccc|ccccccc}
0 & 1 & & & & & & 0 & 1 & & & & \\
0 & \times & 1 & & & & & 0 & \times & 1 & & & \\
0 & \times & \times & 1 & & & & 0 & \times & 1 & 0 & & \\
0 & \times & \times & \times & 1 & & & 0 & \times & \times & \times & 1 & \\
0 & \times & \times & \times & \times & 1 & & 0 & \times & \times & \times & 1 & 0 \\
0 & \times & \times & \times & \times & \times & 1 & 0 & \times & \times & \times & \times & 1 & 0 \\
\end{array}
\right].
\]
This leads to a polynomial of degree $d=30$. Using the locally convergent iterative method described above, we successfully computed a real solution vector. 

For the actual numbers, we refer to the software. 
The maximum absolute value in the coefficient vector for the complex case is $3.7\cdot10^{5}$, and the Jacobian evaluated at the solution has a condition number of $5.8\cdot10^{11}$. The maximum absolute value in the coefficient vector for the real case is $5.0\cdot10^{4}$, and the condition number of the Jacobian at the solution is $2.10\cdot10^{10}$.

The highest degree polynomial set can indeed be different for the real and complex cases, and the simulations suggest the following. 

\begin{conjecture}[Six multiplications] For complex coefficients,
\[
    \max\{d: \Pi_d \subset {\pistarbar{6}}\}=32.
    \]
    For real coefficients,
\[
    \max\{d: \Pi_d \subset {\pistarbar{6}}\}=30.
    \]    
\end{conjecture}

\subsection{Seven multiplications} 
For seven multiplications, we believe that the state of the art is 
$d=30$ \cite[Table~3]{saib21} with $s=4$, which in our context corresponds to the reduction $a_{23}=a_{33}=a_{34}=a_{43}=a_{44}=a_{45}=0$ and $a_{78}=0$. In Figure~\ref{fig:reduced} we see that the highest admissible degree appears to be $d=42$ with $r=5$. Using the same technique as for $m=6$, 
we were able to find a solution to the system for the Taylor expansion of the exponential. We used scaling $\alpha=16$,  complex arithmetic, and the following structure:
\begin{equation}    \label{eq:seven_mult_struct}
[A|B] = \left[
\begin{array}{cccccccc|ccccccccc}
0 & 1 &   &   &   &   &   &   & 0 & 1 &   &   &  \\
0 & \times & 1 &   &   &   &   &   & 0 & \times & 1 &   &  \\
0 & \times & 0 & 1 &   &   &   &   & 0 & \times & 1 & 0 &  \\
0 & \times & \times & \times & 1 &   &   &   & 0 & \times & \times & 1 & 0 &  \\
0 & \times & 0 & \times & \times & 1 &   &   & 0 & \times & \times & \times & 1 & 0 & \\
0 & \times & \times & \times & \times & \times & 1 &   & 0 & \times & \times & \times & \times & 1 & 0 & \\
0 & \times & \times & \times & \times & \times & \times & 1 & 0 & \times & \times & \times & \times & \times & 1 & 0 &  \\
\end{array}
\right].
\end{equation}
The fixed elements $a_{3,3}=a_{5,3}=0$ were selected empirically 
in an attempt to avoid very large condition numbers.
The maximum absolute value in the coefficient vector is $2.8\cdot10^{6}$, and the condition number of the Jacobian at the solution is $2.9\cdot10^{13}$. From this we conclude the following in a complex sense.
\begin{conjecture}[Seven multiplications]
\[
    \max\{d: \Pi_d \subset {\pistarbar{7}}\}=42.
    \]
\end{conjecture}

\begin{remark}[Computational example]\rm
Although the main contributions of this manuscript are theoretical in character, we also wish to point out the practical value of the simulations. For example, the output of the simulations concerning $m=7$ multiplications can be directly used to compute the matrix exponential, in a very competitive way. In fact, in theory (in the sense of number of floating point operations) the method is the fastest for large norm matrices, as far as we know. 
Although a complete computational study is beyond the scope of this paper, we provide (for a specific but random matrix) a comparison with the scaling-and-squaring algorithm \cite{Moler:2003:NINETEEN}, which is the most commonly  method  used in mathematical software for the matrix exponential. In the following we see that our approach is faster or more accurate, or both depending on viewpoint, than two valid
parameter choices for the scaling-and-squaring method. What is marked as \emph{our method} is the evaluation 
\eqref{eq:seven_mult_struct} with coefficients precomputed with the Tikhonov-Newton method.

\begin{verbatim}
julia> Random.seed!(0); setprecision(128); n=200; alpha=16;
julia> A=randn(Complex{BigFloat},n,n); A=20*A/norm(A);
julia> E1 = @btime exp16_deg42_bigfloat(A/alpha); # Our method
  58.648 s (719681983 allocations: 29.50 GiB)  
julia> E2 = @btime exp_sas_6mult_1div(A); # standard version 1 
  65.340 s (762205498 allocations: 31.24 GiB)
julia> E3 = @btime exp_sas_7mult_1div(A); # standard version 2
  80.014 s (860125595 allocations: 35.25 GiB)
julia> # Compute a reference solution with high precision
julia> expAref=mapreduce(i-> (A^i)/factorial(big(i)), +, 0:100);
julia> norm(expAref-E1); # Error for our method
9.25199599560488132729984839075891589943e-33
julia> norm(expAref-E2); # Error for standard implementation version 1
1.935440065066927257455720406200123063236e-30
julia> norm(expAref-E3); # Error for standard implementation version 2
1.14214602689126618911305625822908609328e-36
\end{verbatim}
\end{remark}

\section{Conclusions}\label{sec:conclusions}
This work focuses on a characterization of the  set $\pistar{m}$, with particular attention given to 
minimality 
and to computing the maximum degree polynomial subset of $\pistar{m}$ in the sense of \eqref{eq:maxproblem}.
From our perspective, the minimality question is well understood in this paper.
The determination of the maximum degree polynomial subset can be further investigated. We have only described the case $m\leq 7$ using computational reasoning. Both a theoretical description of the general case, e.g., using further tools from algebraic geometry \cite{Sturmfels:2021:INVITATION}, and a more general computational approach could be of interest and useful in practice.

Based on our simulations, one major component is missing before this can be directly used in matrix function evaluation software: understanding the effect of rounding errors. Evaluating high-degree polynomials is, in general, prone to rounding errors---unless special representations such as a Chebyshev basis are used. In this case, the issue appears even more intricate.
For example, by using high-precision arithmetic, the coefficients were computed such that we could guarantee correctness in full double precision.
However, the fact that the system has a rather large condition number (at least for $m=7$) is an indication that this evaluation is sensitive with respect to these coefficients. Heuristics similar to \cite{jarlebring2024polynomialapproximationsmatrixlogarithm} might be applicable in a general setting. Further work is needed to determine which evaluation schemes, in the continuum of $(A,B,c)$, lead to better numerical stability; a necessary condition for numerical stability is that the condition number is not too large.

The Paterson--Stockmeyer method has proven beneficial not only for evaluating matrix polynomials; similar computational challenges arise in other contexts, such as when the input $X$ is a polynomial. In these scenarios, multiplying two quantities is significantly more computationally demanding than forming linear combinations. The construction in the Paterson--Stockmeyer approach resembles the baby-step giant-step (BSGS) technique introduced in \cite{Shanks:1971:BSGS}, which has found various applications and has been combined with the Paterson--Stockmeyer approach in public key and privacy-preserving cryptography \cite{Han:2020:BOOTSTRAPPING}. Moreover, both the Paterson--Stockmeyer method and BSGS serve as valuable tools in high-precision arithmetic \cite{Johansson:2014:HOLOMONIC,Smith:1989:MULTIPREC}. Open research questions include how the approach presented in this paper, or methods for $\pistar{m}$ in general, can be applied in these contexts.

The fixed-cost computation approach presented in this paper can be complemented by insights from research on composite polynomials or deep polynomials. See \cite{ritt1922prime} for composite polynomials. Rational approximations corresponding to this concept, such as those in \cite{GAWLIK2021105577}, illustrate how successive compositions, e.g.,  $p(f(g(x)))$, can achieve rapid convergence in terms of both the number of compositions and the parameters involved. Similar findings are noted in \cite{wei2025deepunivariatepolynomialconformal},
motivated by the link between this approach and universal approximation in deep learning.
The composite polynomial method can fit within the framework of this paper by zeroing certain elements in matrices $A$ and $B$. Nevertheless, there is a significant distinction in research objectives:
our objective is to minimize the number of matrix-matrix multiplications,
while \cite{wei2025deepunivariatepolynomialconformal} focuses on reducing
the number of parameters. Although some results, such as the approximation of the $p$th root \cite{GAWLIK2021105577}, may be directly applicable, further research is needed to fully explore the differences resulting from these objectives.

\section*{Acknowledgements}
The authors wish to express gratitude for the insightful discussions and feedback from Prof. Kathlén Kohn (KTH Royal Insititute of Technology), particularly regarding Section~\ref{sec:variety}. This research was partially conducted during the first author's sabbatical at EPFL / University of Geneva. The support of the hosts, Prof. Daniel Kressner and Prof. Bart Vandereycken, is greatly appreciated. The authors also greatly acknowledge the valuable comments from Prof. Massimiliano Fasi (University of Leeds) and Prof. Jorge Sastre (Polytechnic University of Valencia).

\bibliographystyle{plain}
\bibliography{eliasbib,misc}

\appendix
\section{Auxiliary material for the proof of Theorem~\ref{thm:dim_2m_conj}} \label{sec:appendix_proof}
The proof is based on two technical lemmas. We note that several of these results hold for arbitrary choices of $(A,B,c)$, but some parts of the statements and some of the proofs are clearer when we assume the specific structure of the triplet $(A,B,c)$ given in \eqref{eq:proof_example}, which is sufficient for the proof of Theorem~\ref{thm:dim_2m_conj}. Note that for our example, $Q_j=A_jB_j$, $A_j=Q_2+Q_{j-1}$, and $B_j=Q_2+\cdots+Q_{j-1}$.

The first lemma states a formula for the degree when we apply a differentiation operator consisting  of a linear combination of derivatives of the elements in $(A,B,c)$ but not of the last rows of $A$ and $B$.
\begin{lemma}\label{lem:recursion}
Let $(A,B,c)$ be of the specific structure given in \eqref{eq:proof_example}.
Let $\DDD_i$ be an operator consisting of a linear combination
of the elements of $A,B$ in rows $1,2,\ldots,i$. Assume that
\[
\deg(\DDD_i Q_3)\le \cdots \le \deg(\DDD_i Q_{i+1}) <\deg(\DDD_i Q_{i+2})
\]
and that $\deg(\DDD_i Q_{i+2}) \ge 1$.
Then, $\deg(\DDD_i Q_{i+2}) < \cdots < \deg(\DDD_i Q_{m+2})$ and 
\begin{equation}\label{eq:recursion_thm_degs}
\deg(\DDD_i p)=\deg(\DDD_i Q_{m+2})=2^{m-1}+\cdots+2^i + \deg(\DDD_i Q_{i+2}).
\end{equation}
\end{lemma}
\begin{proof}
We prove the statement by induction, starting with row $i$. By applying the product rule and using $\DDD_i Q_2=0$, we obtain:

\begin{multline}    
\label{eq:recursion_proof_step1}
\DDD_i Q_{i+3} = (\DDD_i A_{i+3}) B_{i+3} + A_{i+3} (\DDD_i B_{i+3}) = \\
(\DDD_i Q_{i+2}) C_{i+3} + A_{i+3} \left( \DDD_i Q_{3} + \cdots + \DDD_i Q_{i+1} \right)
\end{multline}
where $C_{i+3}=A_{i+3}+B_{i+3}$.
Note that $\deg(C_{i+3}) = \deg(A_{i+3}) = 2^i$. Therefore, on the right-hand side of equation \eqref{eq:recursion_proof_step1}, the degree of the first term is larger than that of the second term, given our assumption that $\deg(\DDD_i Q_{3}) \le \cdots \le \deg(\DDD_i Q_{i+1})  < \deg(\DDD_i Q_{i+2})$. Thus, we have:
\begin{equation}
\label{eq:recursion_proof_deg1}
\deg(\DDD_i Q_{i+3}) = \deg(A_{i+3}) + \deg(\DDD_i Q_{i+2}) = 2^{i} + \deg(\DDD_i Q_{i+2}) > \deg(\DDD_i Q_{i+2}).
\end{equation}
This establishes the base step of our induction.

To proceed with the induction step, assume that the inequality holds for $j$ steps, i.e., 
\[
\deg(\DDD_i Q_{i+2}) < \cdots < \deg(\DDD_i Q_{i+j+1}).
\]
Analogous to the derivation in equation \eqref{eq:recursion_proof_step1}, we have:
\[
\DDD_i Q_{i+j+2} = (\DDD_i Q_{i+j+1}) C_{i+j+2} + A_{i+j+2} \left( \DDD_i Q_{3} + \cdots + \DDD_i Q_{i+j} \right).
\]
Since $C_{i+j+2}$ and $A_{i+j+2}$ have the same degree, and by our induction assumption, the first term has a higher degree, we conclude that:
\begin{equation}
\label{eq:recursion_proof_deg2}    
\deg(\DDD_i Q_{i+j+2}) = \deg(C_{i+j+2}) + \deg(\DDD_i Q_{i+j+1}) = 2^{i+j-1} + \deg(\DDD_i Q_{i+j+1}).
\end{equation}
This completes the proof of the increasing degree progression. The relation \eqref{eq:recursion_thm_degs} follows by applying equation  \eqref{eq:recursion_proof_deg2} for $j = m-i, \ldots, 2$  and using \eqref{eq:recursion_proof_deg1} once.
\end{proof}

The previous lemma (Lemma~\ref{lem:recursion}) gives us the degree of the partial derivative of the output, given the degree of the partial derivative of $Q_{i+2}$. It remains to determine the derivatives of $Q_{i+2}$. We select the differentiation operator in several ways for the purpose of later combining them to form distinct degrees of the Jacobian.
\begin{lemma}\label{lemma:degree-recursion}
Let $(A,B,c)$ be of the specific structure given in \eqref{eq:proof_example}. Then, for $i=2,\ldots,m$ and $j=2,\ldots,i$
\begin{equation}    \label{eq:case2_local}
\frac{\partial Q_{i+2}}{\partial a_{i,j}}
=
Q_j B_{i+2}
\end{equation}
and for $i=3,\ldots,m$ and $j=2,\ldots, i$,
\begin{equation}    \label{eq:case3_local}
    \left( \frac{\partial}{\partial b_{i,j}}-\frac{\partial }{\partial a_{i,j}}\right)Q_{i+2}
    =
    Q_j(A_{i+2} - B_{i+2})
\end{equation}
Moreover, for $i=5,\ldots m$,
\begin{equation}\label{eq:case4a_local}
    \left(
    \frac{\partial }{\partial b_{i-1,i-1}}-
    \frac{\partial }{\partial a_{i-1,i-1}}+2\left(
    \frac{\partial }{\partial b_{i-1,2}}-
    \frac{\partial }{\partial a_{i-1,2}}\right)+
    2\frac{\partial }{\partial a_{i,i}}\right)Q_{i+2}=B_{i+2}B_{i-1}Q_2 + \ordo(X^r),
\end{equation}
where $r={2^{i-1}}$ and also, 
\begin{equation}\label{eq:case4b_local}
    \left(
    \frac{\partial }{\partial b_{m,m}}-
    \frac{\partial }{\partial a_{m,m}}+2\left(
    \frac{\partial }{\partial b_{m,2}}-
    \frac{\partial }{\partial a_{m,2}}\right)+
    \frac{\partial }{\partial c_{m+1}}\right)p=Q_2(2Q_2 - B_{m}).
\end{equation}

\end{lemma}

\begin{proof}
We prove \eqref{eq:case2_local} by applying the product rule: 
\begin{equation}\label{eq:Dqa_local}
    \frac{\partial Q_{i+2}}{\partial a_{i,j}}
    =  \frac{\partial A_{i+2}}{\partial a_{i,j}} B_{i+2}
    + 
    A_{i+2} \frac{\partial B_{i+2}}{\partial a_{i,j}}
    =
    Q_j B_{i+2}.
\end{equation}
Similarly,
\begin{align}\label{eq:Dqb_local}
    \frac{\partial Q_{i+2}}{\partial b_{i,j}}
    =
    Q_j A_{i+2}.
\end{align}
Equation \eqref{eq:case3_local} is an immediate consequence of \eqref{eq:Dqa_local} and \eqref{eq:Dqb_local}.

For notational convenience, we define the differential operator:
\begin{equation}
\DDD_k := 
\frac{\partial }{\partial b_{k,k}}-
\frac{\partial }{\partial a_{k,k}}+2\left(
\frac{\partial }{\partial b_{k,2}}-
\frac{\partial }{\partial a_{k,2}}\right),
\end{equation}
for which we derive the auxiliary relation:
\begin{subequations} \label{eq:case4Dk}
    \begin{align}
    \DDD_kQ_{k+2}
    &=
    (A_{k+2} - B_{k+2})(2Q_{2} + Q_{k}) \\
    &=
    (Q_2 - B_{k+1})  (Q_{2} + A_{k+1}) \\
    &=
    Q_2^2 + Q_2(A_{k+1} - B_{k+1}) - A_{k+1} B_{k+1} \\
    &=2Q_2^2 - Q_2B_{k} - Q_{k+1}.
    \end{align}
\end{subequations}
Here, we have used \eqref{eq:case3_local} with $i=k$, $j=2$ and $j=i$, in the first equality, and the fact that for the specific structure given in \eqref{eq:proof_example}, we have $A_{k+2} - B_{k+2} = Q_2 - B_{k+1}$.

To prove \eqref{eq:case4b_local}, we use \eqref{eq:case4Dk} with $k=m$:
\begin{equation}\label{eq:case4bexpression}
   \left(\DDD_m +\frac{\partial }{\partial c_{m+1}}\right)p=
    2Q_2^2 - Q_2B_{m} - Q_{m+1}+Q_{m+1}=2Q_2^2 - Q_2B_{m}.
\end{equation}

To prove \eqref{eq:case4a_local}, we use that $\DDD_{i-1} Q_2 = \cdots = \DDD_{i-1} Q_{i} = 0$ since these elements are independent of row $i-1$. This, together with the chain rule, yields:

\begin{subequations}\label{eq:local:Diffop1}
\begin{align}
    \DDD_{i-1}Q_{i+2} &= B_{i+2} \DDD_{i-1} A_{i+2} + A_{i+2} \DDD_{i-1} B_{i+2} \\
    &= (A_{i+2} + B_{i+2}) \DDD_{i-1} Q_{i+1} \\
    &= (A_{i+2} - B_{i+2}) \DDD_{i-1} Q_{i+1} + 2 B_{i+2} \DDD_{i-1} Q_{i+1}.
\end{align}
\end{subequations}
The final equality is a reformulation which simplifies the following derivation,
\begin{subequations} \label{eq:localproof1}
\begin{align} 
   \left(\DDD_{i-1} + 2\frac{\partial }{\partial a_{i,i}}\right) Q_{i+2}
    &=
    (A_{i+2} - B_{i+2}) \DDD_{i-1} Q_{i+1} + 2 B_{i+2}\left(\DDD_{i-1} Q_{i+1} + Q_{i} \right) \\
    &=
    (A_{i+2} - B_{i+2}) (2Q_2^2 - Q_2B_{i-1} - Q_{i}) + 2 B_{i+2}\left( 2Q_2^2 + Q_2B_{i-1} \right),
\end{align}
\end{subequations}
where we have used \eqref{eq:case2_local} in the first equality. The degree of the first term in the right hand side is given by
\begin{align}
    \deg((A_{i+2} - B_{i+2}) (2Q_2^2 - Q_2B_{i-1} - Q_{i})) &= \deg(Q_i Q_i) = 2^{i-1},
\end{align}
and the degree of the second term is given by
\begin{align}\label{eq:deg_of_expression_case4a}
    \deg(B_{i+2}B_{i-1}Q_2) = 2^{i-1} + 2^{i-4} + 1.
\end{align}
The theorem conclusion \eqref{eq:case4a_local} follows from Equation \eqref{eq:localproof1} and Equation \eqref{eq:deg_of_expression_case4a}.

\end{proof}

\subsection*{Proof of Theorem~\ref{thm:dim_2m_conj}}
We want to prove the equality in \eqref{eq:dim_2m_eq}; but the upper bound is already given in Corollary~\ref{thm:dim_2m}. To prove that $m^2$ is also a lower bound it is sufficient to find one point, i.e., one triplet, with a Jacobian of rank $m^2$. This is done by finding particular linear combinations of partial derivatives of the entries in $(A,B,c)$ with $m^2$ distinct degrees.
We assume the structure given in \eqref{eq:proof_example}. We construct linear combinations of partial derivatives in 5 different ways. We refer to these linear combinations as Cases 1, 2, 3, 4 and 5.

Case 1: The fact that the evaluation scheme is unreduced means $\deg(Q_{1}) = 1$, and $\deg(Q_{j}) = 2^{j-2}$ for $j=2,\ldots,m+2$. This, together with $\deg\left( \frac{\partial p}{\partial c_j} \right) = \deg(Q_{j})$, directly implies \eqref{eq:deg_ci}.

For Case 2, we consider the operator
\[\DDD_i:=\frac{\partial }{\partial a_{i,j}}, \]
for $i=2,\ldots,m$ and $j=2,\ldots,i$. Since $a_{i,j}$ does not appear in the first $i-1$ rows, we have $\DDD_i Q_2=\cdots =\DDD_i Q_{i+1}=0$. Consequently, the conditions of Lemma~\ref{lem:recursion} are satisfied. Moreover, \eqref{eq:case2_local} implies that
\begin{equation}
\deg(\DDD_i Q_{i+2})=\deg(Q_j B_{i+2}) = \deg(Q_j) + \deg(Q_{i+1}) = 2^{i-1}+2^{j-2}.
\end{equation}
Combining this with \eqref{eq:recursion_thm_degs} yields
\begin{subequations}
    \begin{align}
        \deg(\DDD_i p ) &= \deg(\DDD_i Q_{m+2}) = 2^{m-1} + \cdots + 2^{i} + \deg(\DDD_i Q_{i+2} ) \\
        &= 2^{m-1} + \cdots + 2^{i-1} + 2^{j-2}.
    \end{align}    
\end{subequations}
We conclude \eqref{eq:case2_degformula}.

For Case 3, we consider the operator
\[
\DDD_i:=\left( \frac{\partial}{\partial b_{i,j}} - \frac{\partial}{\partial b_{i,j}} \right),
\]
for $i=3,\ldots,m$ and $j=2,\ldots,i-1$.
The conditions of Lemma~\ref{lem:recursion} are satisfied by the same reasoning as in Case 2.
From \eqref{eq:case3_local} we obtain
\begin{equation}
\deg(\DDD_i Q_{i+2}) = \deg(Q_j (A_{i+2} - B_{i+2})) = \deg(-Q_j Q_i) = 2^{i-2} + 2^{j-2}.
\end{equation}
This, together with \eqref{eq:recursion_thm_degs}, implies \eqref{eq:case3_degformula}.

For Case 4, we consider the operator
\[
\DDD_m:= \left(
    \frac{\partial }{\partial b_{m,m}}-
    \frac{\partial }{\partial a_{m,m}}+2\left(
    \frac{\partial }{\partial b_{m,2}}-
    \frac{\partial }{\partial a_{m,2}}\right)+
    \frac{\partial }{\partial c_{m+1}}\right).
\]
    From equation \eqref{eq:case4a_local} we immediately conclude that
    \begin{align}
        \deg(\DDD_m p) = \deg(\DDD_m Q_{m+2}) = \deg(Q_2 B_m) = 2^{m-3} + 1.
    \end{align}
For Case 5, we consider the operator
\[
\DDD_i:= \left(
    \frac{\partial }{\partial b_{i-1,i-1}}-
    \frac{\partial }{\partial a_{i-1,i-1}}+2\left(
    \frac{\partial }{\partial b_{i-1,2}}-
    \frac{\partial }{\partial a_{i-1,2}}\right)+
    2\frac{\partial }{\partial a_{i,i}}\right),
    \]
for $i=5,\ldots,m$.
We show that this operator satisfies the assumptions of Lemma \ref{lem:recursion}. Firstly, we observe that $\DDD_i Q_2 = \ldots = \DDD_i Q_{i} = 0$ by the same reasoning as in Case 2. Therefore, it is enough to show that
$\deg(\DDD_{i} Q_{i+1}) < \deg(\DDD_{i} Q_{i+2})$.
From \eqref{eq:case4a_local} and \eqref{eq:deg_of_expression_case4a} we have that $\deg(\DDD_i Q_{i+2}) = 2^{i-1} + 2^{i-4} + 1 $. Next, we observe that
\begin{subequations}
        \begin{align}
        \deg \left( \DDD_{i} Q_{i+1} \right)
        &=
        \deg\left(
        \left(
        \frac{\partial }{\partial b_{i-1,i-1}}-
        \frac{\partial }{\partial a_{i-1,i-1}}+2\left(
        \frac{\partial }{\partial b_{i-1,2}}-
        \frac{\partial }{\partial a_{i-1,2}}\right)+
        2\frac{\partial }{\partial a_{i,i}}\right)
        Q_{i+1}
        \right) \\
        &=
        \deg\left((A_{i+1} - B_{i+1})(Q_2 + Q_{i-1}) \right) \\
        &=
        \deg(Q_{i-1}^2) = 2^{i-1} < 2^{i-1} + 2^{i-4} + 1.
    \end{align}
\end{subequations}
Thus the assumptions of Lemma \ref{lem:recursion} are satisfied. It follows that
\begin{subequations}
    \begin{align}
        \DDD_i p  =  \DDD_i Q_{i+2} &= 2^{m-1} + \cdots + 2^{i-1} + 2^{i-4} + 1, \quad i = 5,\ldots,m, \\
        &= 2^{m-1} + \cdots + 2^{i} + 2^{i-3} + 1, \quad i = 4,\ldots,m-1.
    \end{align}
\end{subequations}

To summarize, we have obtained the following degrees through linear combinations of partial derivatives.
\begin{itemize}
    \item Case 1:\; $0,2^{0},2^{1},\ldots,2^{m}$.
    \item Case 2:\; $2^{m-1} + \cdots + 2^i + 2^{i-1} + 2^{j-2},\;\;\textrm{ for }j=2,\ldots,i,\;\;i=2,\ldots,m$.
    \item Case 3:\; $2^{m-1}+\cdots+2^i+2^{i-2}+2^{j-2},\;\;\;\;\;\textrm{ for }j=2,\ldots,i-1,\;\;i=3,\ldots,m$.
    \item Case 4: $2^{m-3} + 1$.
    \item Case 5: $2^{m-1} + \cdots + 2^{i} + 2^{i-3} + 1,\;\;\;\textrm{ for }i = 4,\ldots,m-1.$
\end{itemize}
These form distinct degrees, in the following way. Firstly, we note that when the degrees are expressed as binary numbers, Case 1 and Case 4 always involve one and two nonzeros respectively, making these degrees distinct.
Similarly, when we express the degrees of Case 2, 3 and 5 as binary numbers, they always involve three or more nonzeros, making them distinct from Case 1 and Case 4. For Case 2, Case 3 or Case 5 to coincide, the degrees as binary numbers must have the same number on nonzeros. This corresponds to choosing the same $i$ in each of the formulas, which leads to different degrees. Therefore we conclude that all the above degrees are distinct from each other.

Counting the number of unique degrees, we get: Case 1: $m+1$; Case 2: $m(m-1)/2$; Case 3: $(m-1)(m-2)/2$; Case 4: $1$; Case 5: $m-4$. In total, we have $m^2$ distinct degrees for $m\geq 4$. 

\end{document}

%% file: intro.tex
The application that motivates the research question in this paper, is the computation of matrix functions in the sense of Higham \cite{Higham:2008:MATFUN}, which is a classical problem in numerical linear algebra. 
We define a matrix function as an extension of a scalar function from $f:\CC \rightarrow \CC$ to matrices, i.e., $f:\CC^{n\times n} \rightarrow \CC^{n\times n}$. 
Important matrix functions like $e^X$, $\text{sign}(X)$, and $\sqrt{X}$ are crucial in various contexts such as linear ODEs \cite{Moler:2003:NINETEEN}, control theory \cite{Byers:1987:RICCATI}, network analysis \cite{Estrada:2010:ESTRADA}, and quantum chemistry \cite{Rubensson:2008:DENSITY}; see \cite[Chapter~2]{Higham:2008:MATFUN} for more applications.
Further applications relevant for our setting appear in problems where the action of a matrix function on a vector $b\in\CC^{n}$, i.e., $f(X)b$ is needed, see \cite{Fasi:2024:ChallengesCompMatFun}.

In this paper we study methods for computing matrix functions when the input matrix $X\in\CC^{n\times n}$ is very large and dense. In particular, we consider a family of methods that only utilize the following two operation types:
\begin{itemize}
    \item Linear combination of two matrices $Z\leftarrow \alpha X+\beta Y$
    \item Multiplication of two matrices $Z \leftarrow X\cdot Y$.
\end{itemize}
A direct consequence of considering only these basic operation types is that this family of methods computes matrix polynomials. Another direct consequence is that the first operation type can be viewed as free in terms of computational cost, since the computational complexity is $O(n^2)$ and $O(n^3)$ respectively for the two considered operation types.

We want to study the polynomials that can be computed with a given cost. Since the cost is essentially given by the number of matrix-matrix multiplications, 
methods with a given \emph{fixed cost} corresponding to $m$ multiplications form evaluations of  polynomials in the set
\begin{equation}\label{eq:pistar_org}
    \pistar{m} := \{ p \in \pigrad{2^m}\hspace{-3pt}:\hspace{-1pt}\text{$p(X)$ is computable with $m$ matrix-matrix multiplications} \},
\end{equation}
where $\pideg{d}$ is the set of all polynomials of degree $d$ in the usual sense.
Here, $\pigrad{d}$ is the closure of $\pideg{d}$, that is, the set of 
polynomials of degree $\leq d$.
Although our application stems from matrix polynomials, the set $\pistar{m}$ is, in fact, a univariate polynomial set, not a matrix polynomial set; it is a univariate semi-algebraic set, as we concretize in Section~\ref{sec:eval}.  The set can be defined more abstractly as follows.  Let $\Pi = \KK[x]$ be the polynomial ring over a field $\KK$; typically the polynomial coefficients are scalars with $\KK=\CC$ or $\KK=\RR$, such that  $\Pi_d \subset \Pi$ is the vector space of polynomials in $\KK[x]$ of degree at most $d$. We make the following assumptions about computations involving elements of $\Pi$:

\begin{itemize}
\item Linear combination of two polynomials $r\leftarrow \alpha p+\beta q$
is considered computationally free, where $p,q\in\Pi$ are  polynomials,  and $\alpha, \beta\in\KK$ are scalars.
\item Multiplication of two polynomials $r\leftarrow p\cdot q$, where $p,q\in \Pi$, incurs unit computational cost, independent of the degrees of $p$ and $q$.
\end{itemize}

Now, if we fix a computational budget to $m$, i.e., fix the number of non-scalar multiplications, the total space of computable polynomials is restricted. Since each multiplication can at most double the degree (i.e., $\deg(fg) \le \deg(f) + \deg(g)$), the maximal degree reachable using $m$ multiplications is bounded by $2^m$. The definition of $\Pi^*_{2^m} \subset \Pi_{2^m}$  in \eqref{eq:pistar_org} corresponds to all polynomials in $\Pi$ that can be computed using at most $m$ non-scalar multiplications, and an arbitrary number of free linear combination operations.  Although the matrix polynomial evaluation application is our main source of interest, there are other applications, for example, when the matrix $X$ is replaced by a (scalar) polynomial; see the discussion and references in Section~\ref{sec:conclusions}.
We also wish to stress the difference of this setting in comparison to the approximation theory assumptions. A common application of matrix polynomials involves the approximation of a non-polynomial function $f$.
In such applications, one often wants to compute a polynomial approximation $p$ of $f$ using a \emph{fixed cost}. Hence,
we want to find the best approximation in $\pistar{m}$. 
This is in complete contrast to the classical problem in (scalar) approximation theory, where one seeks the best approximation of a given \emph{fixed degree}.
Since $\pistar{m}\neq \pigrad{2^m}$ in general, the fixed degree and fixed cost approximation problems are different in character. In order to construct methods for the fixed cost approximation problem we need to understand the set $\pistar{m}$. The objective of this paper is to characterize this set. 

%

Techniques to keep the number of matrix-matrix multiplications low have been studied for decades in the context of matrix functions and matrix polynomials.
In Paterson and Stockmeyer's seminal paper \cite{Paterson:1973:OnTheNumber} an algorithm is presented to compute $p(X)$ for $p\in\pigrad{d}$ in $m = \mathcal{O}(\sqrt{d}$) matrix-matrix multiplications, which was made more precise in \cite{Fasi:2018:OptimalityOfPS}. Since they also show that dim$(\pistar{m}) = \mathcal{O}(m^2)$ = $\mathcal{O}(d)$, the Paterson and Stockmeyer algorithm is optimal in the asymptotic order sense. 

The fact that the Paterson and Stockmeyer algorithm is often improvable
has served as the motivation for a number of recent works.
For example, Sastre \cite{Sastre:2018:SMARTPOLY} showed how to improve 
the Paterson--Stockmeyer method in general;  the result shows how to compute most polynomials of degree $8$ using only three multiplications, and
most polynomials of degree $12$ using only four multiplications (in contrast to the Paterson--Stockmeyer algorithm that handles degrees $6$ and $9$, respectively). 
The work has been expanded in various ways (e.g., \cite{saib21}). This research demonstrates how an approximation of a general function, such as the Taylor expansion of the exponential, can be computed using $m=4$ multiplications. Specifically, the approximation known as 15+ corresponds to a polynomial of degree 16 that matches all Taylor coefficients except for the last, i.e., the coefficient for $X^{16}$. Related concepts for rational functions were examined in \cite{sast12}. Further refinements of the algorithms, especially regarding the scaling-and-squaring method \cite[Section~7]{Moler:2003:NINETEEN}, have been investigated for both the matrix exponential \cite{sid19} and the matrix cosine \cite{sird13a} as well as in \cite{Sastre:2015:SMARTPOLY}.
Ideas for efficient polynomial evaluation combined with rational approximations have been used in \cite{bbc19,Blanes:2025:EFFICIENT}, with particular focus on  preservation of Lie algebra properties for the matrix exponential. 
In general, these types of efficient polynomial methods are (at least in theory) better than the standard implementation of the matrix exponential using a Padé approximant combined with the scaling-and-squaring algorithm \cite{Al-Mohy:2010:SANDS}, although further research is needed to definitively support this claim. 
In the terminology of our framework, considerable parts of the results in the literature discussed above, are examples involving polynomials that are elements of 
$\pistar{m}$. Despite this research attention, 
the understanding of this set is far from complete. 
The basis of our analysis is an evaluation scheme (further explained in Section~\ref{sec:eval}) that gives
a parameterization of $\pistar{m}$ with parameters given by a triplet $(A,B,c)$, where $A$ and $B$ are matrices and $c$ a vector. One triplet $(A,B,c)$ corresponds to one polynomial $p\in\pistar{m}$. There are several ways to
obtain the same polynomial; that is, a different triplet $(\hat{A},\hat{B},\hat{c})$ may lead to the same polynomial $p$.

The first set of new results are equivalence transformations. In Section~\ref{sec:transform}, we present procedures
to transform a triplet $(A,B,c)$ into another triplet $(\hat{A},\hat{B},\hat{c})$ that yields the same polynomial $p\in\pistar{m}$.
Several conclusions can be drawn from the transformations. For example, we prove that the first column of the matrices $A$ and $B$ can be selected as zero without loss of generality. More complex transformations reveal that the $(2,2)$ element of $A$ can also be set to zero. Moreover, we establish a transformation related to the third multiplication, which includes an algebraic condition on the elements of $A$ and $B$. Further analysis shows that the $(3,3)$ element of $A$ can be coupled in a simple way to the $(3,3)$ element of $B$ without loss of generality.

The starting point of Section~\ref{sec:minimality} are conclusions of the equivalence transformations.  Using the transformations we can reduce the number of free parameters that parameterize $\pistar{m}$ and conclude that
    \[
    \dim(\pistar{m})\le m^2 \textrm{ for }m\ge 3.
    \]
This bound is sharper than those presented in, for example, \cite{Paterson:1973:OnTheNumber}.  Moreover, we prove that this is optimal and, to our knowledge, our parameterization therefore forms the first minimal parameterization of unreduced evaluation schemes in $\pistar{m}$.

Further results are presented (in Section~\ref{sec:casestudies}) regarding the problem of finding included polynomial subsets, specifically determining
\begin{equation}\label{eq:included_subset}
\max\{d: \Pi_d \subset \pistarbar{m}\}.
\end{equation}
Using the transformations and tools from computational algebraic geometry, we solve this problem for $m=4$ and provide conjectures supported by strongly indicative computational results in high-precision arithmetic for $m=5$, $m=6$, and $m=7$. Our solutions to \eqref{eq:included_subset} are $12$, $20$, $32$, $30$, and $42$ for $m=4$, $m=5$, $m=6$ (complex arithmetic), $m=6$ (real arithmetic), and $m=7$ (complex arithmetic), respectively.
The numerical experiments are reproducible and provided in a publicly available repository, including generated code for Matlab and Julia.

While the simulations focus on the determination of the dimension of $\pistar{m}$ and the solution to \eqref{eq:included_subset}, other aspects of the findings in this paper may be significant beyond these specific research questions.
The transformations themselves are constructive and may be used to improve  various properties of the evaluation scheme,
such as numerical stability. The study of the polynomial subset in Section~\ref{sec:casestudies} is also constructive.
We provide evaluation schemes for the truncated Taylor expansion of the matrix exponential at specified degrees.
For instance, we achieve the degree-$30$ expansion using $m=6$ multiplications, whereas prior works require $7$ multiplications \cite{sid19}. Moreover, our simulations are applicable to several different functions, and the software we provide can be used to
compute the evaluation scheme coefficients 
for polynomials other than those reported in this paper.

%% file: variety.tex
\subsection{The semi-algebraic set perspective of $\pistar{m}$}\label{sec:variety}
For our study, we need concepts from algebraic geometry. 
Let $\KK$ be the field of the parameters in the evaluation.
Let $\XX$ be the product space associated with the parameters. That is,
$\XX:=\KK^{m(m+3)/2}\times\KK^{m(m+3)/2}\times \KK^{m+2}=\KK^s$,
where $s=m^2+4m+2$ is the number of parameters in the evaluation.

If we denote  the evaluation scheme \eqref{eq:generateB}--\eqref{eq:poly_def}
by the map $\Phi:(A,B,c)\mapsto p\in \pigrad{d}$, we can describe $\pistar{m}$ as its image
\begin{equation}\label{eq:Phi}
\Phi(\XX) = \pistar{m}.    
\end{equation}
By construction, $\Phi$ is algebraic, since it depends polynomially on the parameters.
It follows from the Tarski--Seidenberg theorem \cite[Theorem~8.6.6]{mishra1993algorithmic} that $\pistar{m}$ is a semi-algebraic set, since it is the image of an algebraic map over a vector space.

Example \eqref{eq:epsilonexampleAB}--\eqref{eq:epsilonexampleC} illustrates that the set is not necessarily an
algebraic variety, i.e., it is only semialgebraic.
More precisely, let  $\widebar{\cdot}$ denote topological closure in a Zariski sense. Then, in relation to the example \eqref{eq:epsilonexampleAB}--\eqref{eq:epsilonexampleC} we see that 
$\epsilon=0$ is not a valid choice, and indeed $X^7 \not\in \pistar{3}$. However, since $\epsilon=0$ is a limit point, 
we have $X^7\in \widebar{\pistar{3}}$.
 Therefore  $\pistar{m}$ is not a variety since it does not include all of its limit points.

If $\varphi_0,\ldots,\varphi_d$ denotes a basis of the ambient space $\KK[x]_d$, e.g., as in our simulations a monomial basis, and let $J$ be the Jacobian of the map $\Phi$ with an output expressed in this basis. Following the standard definitions in algebraic geometry \cite{shafarevich1994basic-book1}, the dimension of a semi-algebraic set is given by the rank of $J$ at any non-singular point. Most points in the image of a smooth map are non-singular; therefore for almost all $(A,B,c)$ we have
\[
\dim \pistar{m} = \rank J,
\]
where $J\in\KK^{(d+1)\times s}$, $d=2^m$ and $s$ is the number of parameters in the method. 

The dimension of $\pistar{m}$ is always bounded by the number of free parameters, and in this setting it can be concluded from the size of $J$ that
\begin{equation}\label{eq:dim_s_bound}
\dim \pistar{m} \le s=m^2+4m+2.
\end{equation}
In this paper we improve this bound by reducing the number of free parameters and in Section~\ref{sec:minimality} we subsequently conclude that $m^2$ is the exact dimension, for $m>3$.

We note that the properties of the set $\pistar{m}$ differ from an algebraic perspective when considering $\KK=\CC$ versus $\KK=\RR$. Most results in this paper are applicable to both cases; however, when explicitly discussing the Zariski closure, we will assume $\KK=\CC$ unless stated otherwise. In the simulations presented in Section~\ref{sec:casestudies}, we concentrate on $\KK=\CC$ and provide additional observations for $\KK=\RR$.

%% file: gfx/reduced1.tex
   \pgfplotsset{minor grid style={color=blue}} 
\begin{tikzpicture}

\begin{axis}[
     height=220pt,
    xlabel={$m+\epsilon r$},
    ylabel={degree},
    xmin=1, xmax=8,
    ymin=0, ymax=130,
    xtick={1,2,3,4,5,6,7},
    ytick={0,5,...,130},
    legend pos=north west,
    grid style=dashed,
    grid=major,
]

\input{gfx/reduced_data.tex}

\end{axis}

\end{tikzpicture}

%% file: gfx/reduced2.tex
   \pgfplotsset{minor grid style={color=blue}} 
\begin{tikzpicture}

\begin{axis}[
     width=200pt,
     height=220pt,
    xlabel={$m+\epsilon r$},
    ylabel={degree},
    xmin=4.8, xmax=7.9,
    ymin=16, ymax=48,
    xtick={1,2,3,4,5,6,7,8},
    ytick={0,2,...,54},
    legend pos=north west,
    grid style=dashed,
    grid=major,
 ]
\input{gfx/reduced_data.tex}
\legend{}
\end{axis}

\end{tikzpicture}